\documentclass[10pt,a4paper]{amsart}
\usepackage[margin=16mm]{geometry}
\usepackage[utf8]{inputenc}
\usepackage[T1]{fontenc}
\usepackage{amsmath}
\usepackage{amsfonts}
\usepackage{amssymb}
\usepackage[version=4]{mhchem}
\usepackage{stmaryrd}
\usepackage{bbold}
\usepackage{tikz}
\usetikzlibrary{arrows.meta}
\usetikzlibrary{bending}
\usepackage{epsfig}
\usepackage{amssymb,amsmath,amsthm,ytableau}
\usepackage{amscd,pb-diagram}
\usepackage[all]{xypic}
\usepackage[numbers]{natbib}
\usepackage[colorlinks]{hyperref}
\usepackage{MnSymbol}
\usepackage{amsmath, amsfonts, amssymb, amscd, enumerate}
\usepackage{comment}

\def\a{\alpha}
\def\b{\beta}

\def\e{\varepsilon}

\renewcommand{\L}{\Lambda}

\def\i{\iota}

\def\l{\lambda}
\def\o{\omega}

\def\u{\upsilon}

\def\G{\Gamma}

\chardef\tempcat=\the\catcode`\@
\catcode`\@=11
\def\cyracc{\def\u##1{\if \i##1\accent"24 i
    \else \accent"24 ##1\fi }}
\newfam\cyrfam

\DeclareFontFamily{OT1}{msb}{}{}
\DeclareFontShape{OT1}{msb}{m}{n}
 {  <5> <6> <7> <8> <9> <10> gen * msbm
      <10.95><12><14.4><17.28><20.74><24.88>msbm10}{}
\DeclareMathAlphabet{\bubble}{OT1}{msb}{m}{n}

\def\bR{{\mathbb R}}

\newfont{\goth}{eufm10 scaled \magstep1}

\def\ga{\mathfrak{a}}

\def\gg{\mathfrak{g}}

\def\gk{\mathfrak{k}}
\def\gm{\mathfrak{m}}
\def\gn{\mathfrak{n}}

\def\gq{\mathfrak{q}}

\def\gt{\mathfrak{t}}

\def\gz{\mathfrak{z}}

\def\gsl{\mathfrak{sl}}

\def\Sp#1{{\mathrm{Sp(#1)}}}

\newfont{\mcal}{eusm10 scaled \magstep1}

\def\cn{\mbox{\mcal N}}

\def\Id{\mathrm{Id\;}}

\def\End{\mathrm{End\;}}

\def\vol{\mathrm{vol\;}}

\def\Ad{\mathrm{Ad}}
\def\ad{\mathrm{ad}}

\newtheorem{Prop}{Proposition}
\newtheorem{Cor}{Corollary}
\newtheorem{Lem}{Lemma}
\newtheorem{Def}{Definition}

\def\bt{\begin{theorem}}
\def\et{\end{theorem}}
\def\bp{\begin{Prop}}
\def\ep{\end{Prop}}
\def\bc{\begin{Cor}}
\def\ec{\end{Cor}}
\def\bl{\begin{Lem}}
\def\el{\end{Lem}}
\def\bd{\begin{Def}}
\def\ed{\end{Def}}

\def\n{\nabla}

\def\be{\begin{equation}}
\def\ee{\end{equation}}
\def\arr{\begin{array}{rlll}}
\def\ea{\end{array}}
\def\bea{\begin{eqnarray}}
\def\eea{\end{eqnarray}}
\def\bean{\begin{eqnarray*}}
\def\eean{\end{eqnarray*}}

\newtheorem{theorem}{Theorem}[section]

\newtheorem{corollary}{Corollary}[section]

\newtheorem{example}{Example}[section]

\newtheorem{lemma}{Lemma}[section]

\newtheorem{proposition}{Proposition}[section]

\theoremstyle{definition}
\newtheorem{definition}{Definition}[section]
\newtheorem{remark}{Remark}[section]

\newcommand{\Z}{\mathbb{Z}}
 \newcommand{\R}{\mathbb{R}}
\newcommand{\g}{\mathfrak g}
 \newcommand{\gh}{\mathfrak h}
\newcommand{\gl}{\mathfrak{l}}

\newcommand{\Span}[1]{\mathrm{span}\left(#1\right)}
\DeclareMathOperator{\LL}{LGr}

\DeclareMathOperator{\id}{id}

\DeclareMathOperator{\Lie}{Lie}

\DeclareMathOperator{\GL}{\mathsf{GL}}

\DeclareMathOperator{\SL}{\mathsf{SL}}
\DeclareMathOperator{\Gr}{Gr}

\DeclareMathOperator{\Hom}{Hom}
\DeclareMathOperator{\Hess}{Hess}

\DeclareMathOperator{\Smbl}{Smbl}

\DeclareMathOperator{\dd}{d}

\DeclareMathOperator{\Eigen}{\mathcal{A}}
\newcommand{\E}{\mathcal{E}}

\newcommand{\CC}{\mathcal{C}}

\newcommand{\sll}{\mathfrak{sl}}

\newcommand{\Th}{^\textrm{th}}

 \newcommand{\df}{\stackrel{\textrm{def.}}{=}}

\newcommand{\CSUB}{\mathfrak{a}}

\newcommand{\Mcont}{N}
\begin{document}

\title[Bi--Lagrangian effective forms]{Invariant Monge--Amp\`ere equations on contactified para--K\"ahler manifolds}

\author{Dmitri~Alekseevsky}
   \address{Institute for Information Transmission Problems, B. Karetny
per. 19, 127051, Moscow (Russia) and University of Hradec Kralove, Rokitanskeho 62,
Hradec Kralove 50003 (Czech Republic).}
\email{dalekseevsky@iitp.ru}

\author{Gianni Manno}
   \address{Dipartimento di Matematica ``G. L. Lagrange'', Politecnico di Torino, Corso Duca degli Abruzzi, 24, 10129 Torino, Italy.}
    \email{giovanni.manno@polito.it}

\author{Giovanni Moreno}
 \address{Department of Mathematical Methods in Physics,
 Faculty of Physics, University of Warsaw,
ul. Pasteura 5, 02-093 Warszawa, Poland}
 \email{giovanni.moreno@fuw.edu.pl}
\date{\today}


\maketitle

\begin{abstract}
We develop a method for describing invariant  Monge--Amp\`ere equations in the sense of V.~Lychagin and T.~Morimoto (MAE)   on a homogeneous contact manifold  $N$ of a semisimple Lie group  $G$, which is the  \textit{contactification}    of  the    homogeneous  symplectic   manifold
     $M = G/H = \Ad_G Z  \subset \gg$, where $M$ is  the  adjoint orbit of a    splittable closed element  $Z $  of the Lie  algebra  $\gg = \Lie(G)$. The method is then applied to a ten--dimensional semisimple orbit $M$ of the exceptional Lie group $\mathsf{G}_2$ and a complete list of mutually non--equivalent MAEs on $N$ is obtained.
\end{abstract}

\setcounter{tocdepth}{1}
\tableofcontents

 \textbf{Keywords:} 

%

\section*{Introduction}

Monge--Amp\`ere equations  (MAEs) form a  distinguished  class of  nonlinear   second--order PDEs. They  were  introduces by G.~Monge  in 1784 in  his  pioneering   study of  optimal transportation problem  and continued   by A.M.~Amp\`ere  in 1820.\par
The  classical  MAE   has the  form
\begin{equation*}
    \det\Hess u(x) =  f(x, u(x)))\, ,\quad   x \in \R^n\, ,
\end{equation*}
where $\Hess u(x)= D^2u$ is the Hessian  and    $f(,x,u)$ is a given   function.\par 
Numerous applications     in     differential geometry, meteorology, cosmology, hydrodynamics,
 economics,  optimal mass transportation problem,   etc.,   lead to consideration of a more general class of  MAEs,  given by
\begin{equation}\label{eqEQ1Dmitri}
  \det[\Hess u-A(x, u, Du)] = f(x, u, Du),
\end{equation}  
where  $A $  is a symmetric matrix. \par       
 Monge--Amp\`ere  equations   are intensively studied \cite{Trudinger2008TheME,Alessio2014,Villani2009}.  Many  deep results    about  existence   and unicity of    solutions    are obtain   in the case  when     $\Hess u(x)$ is positively defined,  and as such it  may be  considered as a Riemannian metric  (the so called Hessian metric), see \cite{Mealy1991,Reese2017}.\par 

A  short  history of MAEs     and     their    complex  and  quaternionic    versions  can be found  in the   notes~\cite{Verbitsky2010}  of M.~Verbitsky.\par

V.V. Lychagin \cite{Lyc78,Lychagin1979}  and T.  Morimoto \cite{Morimoto1979} proposed a    general construction of a   wide  class of  MAEs     in terms   of     contact    and   symplectic  geometry.
Let  $(N,D)$ be  a $2n+1$ dimensional  manifold with contact structure  $D$, i.e., a  codimension--one  distribution locally defined  by 1--form $\theta$  with   $d\theta^n \wedge \theta  \neq 0$.  Then  any $n$--form
    $\Omega$ defines a MAE  $  \E_{\Omega}$.   A  solution of  the  equation is a Legendrian  submanifold  $L \subset N$,   $ \theta|_{L} =0$,  which annihilates $\Omega$.
By    Darboux  theorem,    $(N, \theta )$ is locally identified  with  the  space  $J^1F$  of   1-jets  of functions  on  an $n$--dimensional manifold   $F$.  Then, in  terms    standard local  coordinates  $x^i, p_j, u$, the  MAE   reduces to  the   equation \eqref{eqEQ1Dmitri}  on   a  function    $u(x)$. More precisely,
 1-jet $j^1u(x) \subset N= J^1F$  is    a  Legendrian   submanifold    and  it is  a solution of MAE  if the   function $u(x)$     satisfies  the        equation  \eqref{eqEQ1Dmitri}.\par 
 If  the  $n$--form  $\Omega$  does not depend  of the coordinate $u$, it  can be   considered  as  a   $n$--form on the symplectic manifold $T^*F$  and  then    solutions  of  the   equation  $E_{\Omega} \subset  T^*F$ are  Lagrangian   submanifolds   of the  symplectic manifold  $T^*F$ which  annihilate  $\Omega$.\par 
 Starting  from the paper by G.~Monge, the  most  important  application of MAEs  remains the  application to the  optimal transportation  problem  and   related  problems.\par 
New  geometric approach to this subject had been developed  in  a series of papers \cite{Kim2010,WARREN2010,WARREN2011816,kim2010pseudoriemannian,HarveyBlaine2012} by Y.-H.~Kim, R.~McCann, M.~Warren,  and Harvey-Lawson. They  established  a  closed   relationship between the classical Monge--Kantorovich mass transportation
problem   and   para--K\"ahler  geometry.  More precisely, they proved  that   under some  assumptions    the  solution of a  Monge--Kantorovich   problem   reduces  to the   construction of a special  Lagrangian  submanifold  $L$     in      some  $2n$--dimensional   para--K\"ahler  manifold  $(M,g, \omega)$  with para--holomorphic   $n$--form  $\Phi$,  that is  an  $n$--dimensional  real Lagrangian
$(\omega|_L =0)$   submanifold which annihilate     the  real $n$--form  $\Omega=\mathrm{Im} \Phi$.\par 

In  this  paper   we develop       an  approach  for  describing      invariant Monge--Amp\`ere   equations  in the  sense of Lychagin--Morimoto on  homogeneous   contact  manifolds $N$ of semisimple Lie  group  $G$.  More precisely,  we   consider     manifolds  $N=  G/L$  that  are  contactifications  of    homogeneous para--K\"ahler manifolds   $M =  G/K$,  described  by  \cite{Dmitrii_V_Alekseevsky_2009,Alekseevsky20078}.  The  method  is  applied   to   the   classification of   all  invariant Monge--Amp\`ere  equations  on  the 11--dimensional    contact manifold 
$N= G_2/\SL(2,\R)$  of the  exceptional    non--compact Lie group $\mathsf{G}_2$,    associated to the 10--dimensional  para--K\"ahler  real  flag manifold $\mathsf{G}_2/\GL(2, \R)$.\par 
Note  that     invariant second--order  PDEs    with $\mathsf{G}_2$  symmetries   had been studied    by Yamaguchi  \cite{Yamaguchi1999} and    a    description     of   invariant PDEs    with  $\mathsf{G}_2$
  symmetry, different from  MAE    was   obtained  in   the  remarkable  paper    \cite{The2018} by Dennis~The.  Classifications  of  different   classes of  MAEs      were   given in the  papers
\cite{doi:10.1142/9789814354394_0010,ASENS_1993_4_26_3_281_0,MR2352610,Kruglikov1999,imsatfia2013application} 
and the problem of  equivalence    and  the   conditions of linearizability of MAE   was  treated  in 
\cite{Kushner2008,Kushner2010OnCE}.

\subsection*{Structure of the paper}
In Section~\ref{secBiLagANDParaKahler} we remind that the notion of a bi--Lagrangian manifold is equivalent to the notion of a para--K\"ahler manifold; then we pass to the $G$--homogeneous case and recall some useful results: an adjoint orbit $\Ad_G(Z)$ possesses a bi--Lagrangian structure if and only if the element $Z$ is splittable, and all the bi--Lagrangian structures on $\Ad_G(Z)$ are in one--to--one correspondence with the fundamental $\Z$--gradations of $\gg=\Lie(G)$. Next, we introduce the contactification of a homogeneous para--K\"ahler manifold and, finally, by employing the generalized Gauss decomposition, we prove a theorem that allows to locally identify a para-K\"ahler manifold with the cotangent bundle to a suitable flag manifold.\par
In Section~\ref{secMAEs} we recall the notion of a Monge-Amp\`ere equation in the sense of Lychagin--Morimoto and we extend it to a general contact manifold, thus obtaining the general Monge--Amp\`ere equation associated to an $n$--form (general MAE). The local expression of general MAEs in Darboux coordinates is obtained later, together with the natural interpretation of the fibers of a MAE as hyperplane sections of a Lagrangian Grassmannian.\par
In Section~\ref{secExamplesPK} we work out some examples of fundamental gradations: in particular, we show all fundamental gradations of the algebra $\sll(V)$ and $\gg_2=\Lie(\mathsf{G}_2)$.\par
In Section~\ref{secG2_case} we focus on the 10--dimensional $\mathsf{G}_2$--homogeneous manifold $M=\mathsf{G}_2/\mathsf{GL}_2(\R)$ and we construct a basis of the space of invariant effective 5--forms on the contactification $N$ of $M$. This allows to provide a list of all $\mathsf{G}_2$--invariant MAEs on $N$, that are second--order (nonlinear) PDEs in 5 independent variables; we finally  establish which of them are contact--equivalent.

\subsection*{Acknowledgments}
G.~Manno gratefully acknowledges support by the project  ``Finanziamento alla Ricerca'' under the contract numbers \texttt{53\_RBA17MANGIO}, \texttt{53\_RBA21MANGIO},  and by the
PRIN project 2017 ``Real and Complex Manifolds: Topology, Geometry and
holomorphic dynamics'' (code \texttt{2017JZ2SW5}).  G.~Manno is a member of\linebreak GNSAGA
of INdAM.  G.~Moreno is supported by the Polish National Science Centre    project “Complex contact
manifolds and geometry of secants” \texttt{2017/26/E/ST1/00231}.

  




\section{Homogeneous para--K\"ahler manifolds and their contactification}\label{secBiLagANDParaKahler}

\subsection{Bi--Lagrangian and para-K\"ahler  structures}
Below we introduce two equivalent categories made of the objects we will be working with.
\begin{definition}
    An \textit{almost para--complex structure} on a $2n$--dimensional manifold $M$ is a decomposition 
\begin{equation}\label{eqBiLagStr}
    TM = T^+M \oplus  T^-M
\end{equation}
of the tangent bundle  $TM$ into the direct sum of  two $n$--dimensional distributions $T^{\pm}M$ or, equivalently,   a field $I\in\Gamma(\End(TM))$ of endomorphisms, such that $I|_{T^{\pm}M}=\pm\id_{T^{\pm}M}$. An almost para--complex structure is called a \textit{para--complex structure} if the distributions $T^{\pm}M$ are
integrable, i.e., it holds
\begin{equation}\label{eqIntTpmM}
    [\Gamma(T^{\pm}M),\Gamma(T^{\pm}M)]\subset \Gamma(T^{\pm}M)\, .
\end{equation}
\end{definition}
\begin{definition}\label{defBiLag}
The    decomposition~\eqref{eqBiLagStr}  on    a symplectic manifold   $(M,\omega)$ is called     \textit{almost  bi--Lagrangian}  if  $ \omega|_{T^{\pm }_xM}$ vanishes identically at any point  $x \in M$.
If, moreover,   the  distributions  $T^{\pm}M$ are  integrable, it is called a \textit{bi--Lagrangian structure}.  
\end{definition}
The  integrable  submanifolds (of maximal dimension)   of  the distributions $T^\pm M$ of a bi--Lagrangian structure on a $2n$--dimensional symplectic manifold $(M,\omega)$  are Lagrangian submanifolds of $M$, i.e., they are $n$--dimensional and $\omega$ vanishes identically on them; see also~\cite{Bryant_2001}.

\begin{definition}
An \textit{almost para--Hermitian manifold} is a pseudo--Riemannian manifold $(M, g)$ equipped with an almost para--complex structure $I$, such that the  distributions $T^\pm M$
are  (absolutely) isotropic, i.e.,  $g|_{T_x^{\pm} M}$ vanishes indentically at any point $x \in M$. 
\end{definition}
Given an almost para--Hermitian manifold $(M,g,I)$, the (skew--symmetric)  two--form   
\begin{equation}\label{eqDefOmega}
    \omega: = g(\,  \cdot\,  , I(\, \cdot\, ))\, .
\end{equation}
is called the \textit{K\"ahler  form}.
\begin{definition}\label{defPK}
An almost para--Hermitian manifold $(M,g,I)$, such that   $I$ is  a para--complex     structure,    is called  a  \textit{para--Hermitian manifold}.    A para--Hermitian manifold   $(M, g, I)$ is called   a \textit{para--K\"ahler manifold}  if  the  para--complex   structure     $I$   is parallel with respect to the Levi--Civita connection $\nabla^g$ of $g$, i.e., $\nabla^g I =0$.
\end{definition}




The next well--known results show that the  notion of a bi--Lagrangian manifold  $(M,\omega, I)$  is  equivalent to the notion of a para-K\"ahler manifold   $(M, g, I)$. 
 \begin{lemma}\label{lemKuszul}
 Let $(M, g, I)$ be a para--Hermitian manifold: then  $\nabla^g I =0$ if and only if the  K\"ahler  form  $\omega$ given by~\eqref{eqDefOmega} is closed, i.e.,  $d \omega =0$.
 \end{lemma}
 \begin{proof}
It suffices to show that the identities
    \begin{eqnarray}
        d\omega (v,w,u) &=& g(\nabla^g_w(I)(v),u)+g(\nabla^g_u(I)(w),v)+g(\nabla^g_v(I)(u),w)\, ,\label{eqRompicapo1}\\
        d\omega(v,w,u)+d\omega(v,I(w),I(u))&=&-2g(\nabla^g_v(I)(w),u)\, ,\label{eqRompicapo2}
    \end{eqnarray}
hold for all commuting vector fields $v,w,u,I(w),I(u)$ on $M$, see~\cite[Proposition~4.16]{Ballmann2006}.\par 
In view of such a commutativity and of~\eqref{eqDefOmega}, we obtain
\begin{eqnarray*}
    d\omega (v,w,u)&=&v(\omega(w,u))+w(\omega(u,v))+u(\omega(v,w))\\
&=&v(g(w,I(u)))+w(g(u,I(v)))+u(g(v,I(w)))\, ,\\
&=&g(\nabla^g_v(w),I(u))+g(w, \nabla^g_v(I(u)))+\ldots \\
&=&g(\nabla^g_v(w),I(u))+g(w, \nabla^g_v(I)(u))+g(w, I(\nabla^g_v(u)))+\ldots \\
&=&g(w, \nabla^g_v(I)(u))+\omega(\nabla^g_v(w),u)+\omega(w,\nabla^g_v(u))+\ldots\\
&=&g(\nabla^g_w(I)(v),u)+g(\nabla^g_u(I)(w),v)+g(\nabla^g_v(I)(u),w)\, ,
\end{eqnarray*}
where the dots denote cyclic permutations of $(v,w,u)$, that is formula~\eqref{eqRompicapo1}.\par
Let us observe now that, in view of the Leibniz rule for the covariant derivative and the fundamental identity $\nabla^g(g)=0$, it holds 
\begin{eqnarray}
    2g(\nabla^g_v(I)(w),u)&=& 2g(\nabla^g_v(I(w)),u)-2g(I(\nabla^g_v(w)),u)=\nonumber\\
    &=&2g(\nabla^g_v(I(w)),u)+2g(\nabla^g_v(w),I(u))\, .\label{eqPrimaDiKoszul}
\end{eqnarray}
The Koszul formula, applied to both the addends of~\eqref{eqPrimaDiKoszul}, yields
\begin{eqnarray}
 2 g(\nabla^g_v(I(w)),u)  &=& v(g(I(w),u))+I(w)(g(v,u))-u(g(v,I(w)))=\nonumber\\
    &=&v(\omega(u,w))+I(w)(\omega(v,I(u)))+u(\omega(w,v)))\label{eqKoszul1}
\end{eqnarray}
and
\begin{eqnarray}
    2g(\nabla^g_v(w),I(u))&=& v(g(w,I(u)))+w(g(v,I(u)))-I(u)(g(v,w))\nonumber=\\
    &=&v(\omega(w,u))+w(\omega(v,u))+I(u)(\omega(I(w),v)))\, ,\label{eqKoszul2}
\end{eqnarray}
respectively. By taking the sum of~\eqref{eqKoszul1} and~\eqref{eqKoszul2} we obtain
\begin{eqnarray}
    2g(\nabla^g_v(I)(w),u)&=&v(\omega(u,w))+I(w)(\omega(v,I(u)))+u(\omega(w,v)))+v(\omega(w,u))+w(\omega(v,u))+I(u)(\omega(I(w),v)))=\nonumber\\
    &=& I(w)(\omega(v,I(u)))+u(\omega(w,v))) +w(\omega(v,u))+I(u)(\omega(I(w),v)))\, .\label{eqKoszul3}
\end{eqnarray}
It is then easy to see that
\begin{eqnarray*}
d\omega(v,w,u)+d\omega(v,I(w),I(u))&=&v(\omega(w,u))+w(\omega(u,v))+u(\omega(v,w)))+\\
    &&+v(\omega(I(w),I(u)))+I(w)(\omega(I(u),v))+I(u)(\omega(v,I(w))))=\\
    &=&v(\omega(w,u))+w(\omega(u,v))+u(\omega(v,w))+\\
    &&-v(\omega(w,u))+I(w)(\omega(I(u),v))+I(u)(\omega(v,I(w))))\, ,
\end{eqnarray*}
after switching the arguments in all the instances  of $\omega$, coincides with~\eqref{eqKoszul3}: this proves~\eqref{eqRompicapo2}.
 \end{proof}

\begin{proposition}\label{thEquiCat}
If  $(M,g, I)$ is a para--K\"ahler  manifold, then   $(M,\omega,I) $ is a bi--Lagrangian manifold with $\omega$ given by~\eqref{eqDefOmega}; conversely,    if 
 $(M,\omega,I)$ is a  bi--Lagrangian manifold, then  $(M,g,I)$ is a   para--K\"ahler manifold, with $g$ given by 
 \begin{equation}\label{eqDefGi}
     g:=-\omega\circ I\, .
 \end{equation}
\end{proposition}

\begin{proof}
If a  {para--K\"ahler manifold}   $(M,g, I)$ is given, then the form  $\omega$
defined  by~\eqref{eqDefOmega} is symplectic thanks to Lemma~\ref{lemKuszul}; it only remains to show that $T^\pm M$ is Lagrangian, but this is an immediate consequence of    $T^\pm M$ being   $g$--isotropic:
\begin{equation}\label{eqLagrangianita}
    v,w\in T^\pm M\Rightarrow\omega(v,w)=g(v,I(w))=g(v,\pm w)=0\, .
\end{equation}
If a bi--Lagrangian manifold on $(M, \omega, I)$ is given, then one has to show that~\eqref{eqDefGi} defines indeed a pseudo--Riemanian metric, such that $T^\pm M$ are $g$--isotropic and $I$ is $g$--parallel.\par
The first claim is easy and formally mirrors~\eqref{eqLagrangianita}:
\begin{eqnarray*}
v,w\in T^\pm M&\Rightarrow &g(v,w)=\omega(I(v),I(w))=-\omega(\pm v,\pm w)=-\omega(v,w)=0\, ,\\
v\in T^\pm M\, ,w\in T^\mp M&\Rightarrow &g(v,w)=-\omega(\pm v,\mp w)=\omega(\mp  v,\pm w)=-\omega(\pm w,\mp v)=g(w,v)\, ,
\end{eqnarray*}
because $T^\pm M$ is Lagrangian.\par
The second claim follows again from Lemma~\ref{lemKuszul}.
\end{proof}

\subsection{Classification  of  para-K\"ahler   homogeneous manifolds of   a semisimple Lie  group}
In what follows $G$ will denote a real Lie group and $Z$ an  element of its Lie algebra $\gg:=\Lie(G)$; the Killing form on $\gg$ will be denoted by the symbol $B$.
 

     \begin{definition}
An element $Z \in \gg$ is called  \textit{semisimple} if such is its adjoint operator $\ad_Z$. 
          A semisimple   element $Z \in \gg$ is called  \textit{splittable} if   $\ad_Z$ has   real eigenvalues;   it  is called   \textit{closed}, if it generates   a   closed 1--parameter  subgroup  $\{\exp t Z\mid t\in\R\} \simeq \bR $ of    $G$.     
     \end{definition}

%
It is a classical result by Kirillov, Kostant and Souriou   that,  up to a  central  extension,   any  homogeneous symplectic  manifold   $M = G/H$  is a  coadjoint orbit in the dual $\gg^*$ of the Lie algebra $\gg$ of $G$.  If $G$ is   semisimple, then the Killing  form $B$ is not degenerate and therefore the  coadjoint orbit    $\Ad_G \xi \subset \gg^*$  can be identified with the  adjoint  orbit 
$\Ad_G Z  \subset \gg$, with     $Z$ given by $Z=  B^{-1}\circ \xi$.
    \begin{theorem}[Kirillov--Kostant--Souriou] 
    Up to a covering, any  homogeneous   symplectic    manifold    $M = (G/H, \omega)$ of a semisimple Lie group $G$    is isomorphic   to an adjoint orbit   $\Ad_G Z \subset  \gg$ equipped   with the $G$--invariant symplectic form $\omega$ given by 
\begin{equation}\label{eqKKSform}
    \omega_x(X, Y) := B(x, [\Ad_{x^{-1}\ast}(X), \Ad_{x^{-1}\ast}(Y)])\, ,
\end{equation} 
 for any $x \in M \subset  \gg$ and any $X, Y\in T_xM\subset T_xG$, where $\Ad_{x^{-1}\ast}:T_xG\to T_eG=\gg$ is the tangent map  of $\Ad_{x^{-1}}$ at the point $x$.
\end{theorem}
%
The Kirillov--Konstant--Souriou form~\eqref{eqKKSform} is one of the  many non--equivalent $G$--invariant symplectic  forms on an  adjoint orbit. Let $M=\Ad_GZ=G/H$, where $\gh=\Lie(H)=C_{\gg}(Z)$ is the centralizer of $Z\in\gg$ in the Lie algebra $\gg$ of $G$, and fix a complementary subspace $\gm$ of $C_{\gg}(Z)$ in $\g$:
\begin{equation}\label{eqGiCiGiZeta}
    \gg = C_{\gg}(Z) + \gm \, .
\end{equation}
Then the tangent space $T_Z M$ identifies with $\gm$.
\begin{definition}Let $G$ be a semisimple Lie group  and let $\gm$ be a $B$--orthogonal complement of $C_{\g}(Z)$: then the unique $G$--invariant symplectic structure $\omega$ on $M$, such that
\begin{equation}\label{eqSympFormAdZ}
    \omega_Z (X,Y)   = B(Z, [X,Y])\, ,\quad\forall X,Y \in T_Z M\, ,
\end{equation}
is called \textit{the symplectic structure  on $M=\Ad_GZ$ associated with  $Z\in\gg$}.
\end{definition}

%
\begin{theorem}[Hou-Deng-Kaneyuki-Nishiyama~\cite{Hou1999,Dmitrii_V_Alekseevsky_2009}]\label{thCinesi}
Let $G$ be a semisimple real Lie group and $(M = Ad_GZ, \omega)$ an adjoint orbit of an element $Z \in g = \Lie(G)$, equipped with the invariant symplectic structure $\omega$    associated  with  $Z$. Then the manifold $M$ admits a G--invariant integrable bi--Lagrangian (or, equivalently, para--K\"ahler) structure if and only if $Z$ is a splittable element.
\end{theorem}
Indeed, if $Z $ is splittable, then all the eigenvalues of the operator $\ad_Z $ are real, so that
\begin{equation}
     \gg = \sum_{j\in\Eigen} \gg_j \, ,
\end{equation}
where   $\gg_j$ denotes  the   eigenspace  of  $\ad_Z$ that corresponds to the eigenvalue  $j$ and  $\Eigen:=\{ j\in\R\mid \det(\ad_Z-j\id_\g)=0\}$ is the set of real eigenvalues of $\ad_Z$. In turn, this allows to refine the $B$--orthogonal  reductive decomposition~\eqref{eqGiCiGiZeta} as follows:
\begin{equation}\label{eqGenGaussAlg}
     \gg =\gg^- + \gg^0 + \gn^+  =  \gn^- + \gh + \gn^+ \, ,
\end{equation}
where $\gg^0=\gh =C_\g(Z)$, and
\begin{equation}
     \label{eqDecGiPiuMeno}\gg^\pm=\gn^\pm := \sum_{j\in\Eigen\cap\R^\pm} \gg_j \, .
\end{equation}
The $\ad_\gh$--invariant decomposition~\eqref{eqGenGaussAlg} extends to a    $G$--invariant     decomposition   $TM = T^-M + T^+M $ of $TM$, which turns out to be bi--Lagrangian, because $\gn^\pm$ is a sub--algebra (which ensures integrability) and from the definition~\eqref{eqSympFormAdZ} of $\omega_Z$ it follows that the symplectic form $\omega_Z$ vanishes on $\gn^\pm$.\par 
 
                The $\ad_{\gh}$--invariant decomposition~\eqref{eqGenGaussAlg} is called \emph{a generalized  Gauss   decomposition}, see Definition~\ref{defGaussDec} below.
\begin{corollary}\label{corIntegralsOrbits}
Let $N^\pm$ be the nilpotent subgroup generated by the Lie algebra $\gn^\pm$: the integral submanifold of the distribution $T^\pm M$ passing through the point $x = gZ \in  M$ is given by $gN^\pm Z$.
\end{corollary}
    
         %
     
The  classification  of homogeneous  bi--Lagrangian manifolds of a (complex or real) semisimple Lie group $G$ can be reduced to the description of the fundamental gradations of the corresponding Lie algebra $\gg$. \par 
    A $\Z$--gradation of a Lie algebra $\gg$  is a decomposition   
\begin{equation}
    \gg = \gg^{-k}  + \ldots + \gg^{-1} + \gg^0+ \gg^l + \ldots + \gg^k \,  ,  \label{eqZgrad}
\end{equation} 
satisfying  $[\gg^i,\gg^j]\subset\gg^{i+j}$ for all $i,j\in\Z$.
\begin{definition}
A $\Z$--gradation~\eqref{eqZgrad}  is called \textit{fundamental}  if the subalgebra  
\begin{equation}
    \gg^- :=  \gg^{-k}  + \ldots + \gg^{-1} 
\end{equation} 
is generated  by $\gg^{-1}$.
\end{definition}
Theorem~\ref{thAM1} below    follows  from the  remark that   the operator    $D$  given, for any $j$, by  $D|_{\gg^j} := j \id_{\gg^j}$    is a  derivation of   the Lie  algebra \eqref{eqZgrad},  and  any derivation of  a semisimple Lie  algebra    is inner, i.e., $D = \ad_d$, with $d \in \gg$.\par

\begin{theorem}[Alekseevsky-Medori~\cite{Alekseevsky20078}]\label{thAM1} Let $(M = \Ad_G(Z) = G/H,\omega)$ be as in Theorem~\ref{thCinesi} and let $\gh = C_{\gg} (Z) = \Lie(H)$. There exists a natural one--to--one correspondence between 
\begin{itemize}
    \item[i)] invariant bi-Lagrangian structures $TM = T^+M + T^-M$; 
\item[ii)] $H$--invariant decompositions (called bi--isotropic) of the Lie algebra 
\begin{equation}\label{eqBiIso}
    \gg=\gn^-+\gh+\gn^+ \, , 
\end{equation}
where $\gn^{\pm}$ given by~\eqref{eqDecGiPiuMeno} are subalgebras such that $B|_{\gn^{\pm}}=0$; 
\item[iii)] fundamental $H$--invariant $\Z$--gradations
 \eqref{eqZgrad} with $\gg^0 = \gh$. 
\end{itemize}
\end{theorem}
       More precisely, the bi--isotropic decomposition~\eqref{eqBiIso}, which corresponds to the  fundamental gradation~\eqref{eqZgrad}, is given by 
       \begin{equation}
           \gn^{\pm} = \sum_{\pm i>0}\gg^i
       \end{equation}
       and $\gh  = \gg^0$, while 
the bi--Lagrangian decomposition $TM =T^+M + T^-M$ associated with~\eqref{eqBiIso}  is 
the natural invariant extension of the $H$--invariant decomposition $T_oM =\gn^+ + \gn^-$ of the tangent space of $M = G/H$, $o = eH = [e]_H$, under the standard identification $T_oM = \gg/\gh= \gn^+  +\gn^-$. 

\subsection{Contactification of  a homogeneous para--K\"ahler manifold}\label{secContactification}
    Let $G$ be a real semisimple Lie group    and   $M := \Ad_G Z =G/H$    the  adjoint orbit   of a closed and splittable semisimple element $Z \in \gg$: consider  the   $B$-orthogonal decomposition  $\gh = \gl + \bR Z $ of the stability subalgebra $\gh=\Lie(H)$, where  $\gl=(\bR Z)^\perp$.
 \begin{theorem}
      The  Lie  algebra $\gl$   generates a closed  subgroup $L \subset G$   and the  homogeneous manifold  $\Mcont := G/L $   has the  reductive decomposition
       $$   \gg = \gl + \gn = \gl + (\bR Z + \gm).$$
       The $Ad_L$-invariant 1-form $\theta = B \circ Z$  defines  an invariant  contact structure on $\Mcont$ with the   invariant contact   form $\theta$  which is   an  invariant  extension  of the   1-form $\theta = B \circ Z$.
 \end{theorem}

The natural projection
    \begin{equation}
        \pi:\Mcont = G/L  \longrightarrow M = G/H \, 
    \end{equation}
is a   $G$-equivariant principal bundle with structure group $\bR = \{\exp t Z\mid t\in\R\}$, and $H$ identifies with  the product $H = L\cdot \{\exp t Z\mid t\in\R\}$; moreover,  the contact form $\theta:TN\longrightarrow\R$ turns out to be the connection form of a principal connection on $\pi$, whose curvature form is the symplectic form $\omega=d\theta$ on $M$,  
see~\cite[Theorem~4]{Boothby-Wang_1958}.

\begin{definition}
      The  $G$--homogeneous   contact manifold  $(\Mcont, \theta)$ is called the \emph{contactification} of  the  $G$--homogeneous  para--K\"ahler manifold $(M, \o)$.
\end{definition}


\subsection{Homogeneous  para--K\"ahler manifolds  as a  completion of the  cotangent bundle of  a  real flag manifold}

            Let $(M,\o) $  be   as in Section~\ref{secContactification} above; recall also the bi--isotropic decomposition~\eqref{eqBiIso} from Theorem~\ref{thAM1} and the subgroups $N^{\pm}$ from Corollary~\ref{corIntegralsOrbits}. 
Letting   $P^{\pm} := H \cdot N^{\pm}$ be the closed real  parabolic subgroups of  $G$  with nilradical
             $N^{\pm}$, which is a  semidirect product  $P^{\pm} = H \ltimes N^{\pm}$, we have the following well--known result.
\begin{proposition}\label{propGaussDec}
   The  map
    \begin{eqnarray*}
        N^- \times H \times N^+ &\longrightarrow& G\, ,\\
               (n^-, h, n^+) &\longmapsto& g= n^- h n^+\, , 
    \end{eqnarray*}
             is a diffeomorphism  onto an  open dense    submanifold $\widetilde{G} \subset G$
\end{proposition}
             \begin{definition}\label{defGaussDec}
              The  decomposition
               $$\widetilde{G} =  N^- H N^+$$
                is called the  \textit{generalized Gauss   decomposition}.                \end{definition}
We shall need the real flag manifold $F := G/P^+$ corresponding to the parabolic subgroup $P^+$, together with the   natural projection
\begin{equation}
    \pi : M = G/H \longrightarrow F = G/P^+\, .
\end{equation}
Let $\omega_{\mathrm{std}}$ denote the standard symplectic form on the cotangent space to a smooth manifold.
\begin{theorem}\label{thDarboux}
Up to a zero--measure subset, the   homogeneous  symplectic manifold  $M$ is symplectomorphic to the symplectic manifold $T^*F$, equipped with the  symplectic form $\omega_{T^*F}$ unambiguously defined by
\begin{equation}\label{eqDeformedSymplecticForm}
    \omega_{T^*F}|_{\gg^j\oplus\gg^{-j}}=j\omega_{\mathrm{std}}\, ,\quad\forall j\in\mathbf{N}\, .
\end{equation}
\end{theorem}
\begin{proof}
Let $ o =[e]_H= eH$ be the origin of $M$ and let us take  $\pi(o) =  e P^+ $ as the origin  of $F$: this means that $o\in M=G/H$ corresponds to $Z\in M=\Ad_G(Z)$. The following $n$--dimensional manifolds are then related to each other by an $N^-$--equivariant diffeomorphism:
\begin{eqnarray*}
    N_M&:=&N^-o=\Ad_{N^-}(Z)\subset M\, ,\\
    N_F&:=&N^-\pi(o)= \widetilde{G}/P^+\subset F\, .
\end{eqnarray*}
Indeed, $\pi: M\longrightarrow F$ restricts to an $N^-$--equivariant diffeomorphism $\pi: N_M\longrightarrow N_F$; moreover, by Proposition~\ref{propGaussDec}, the    orbit  $N_F$  is  an   open dense  submanifold   of   $F$.\par 
Next, we prove that the $N^-$--homogeneous vector bundles $T^*N_M$ and $T^+M|_{N_M}$ can be identified by means of the ($G$--invariant) isomorphism between $T^*M$ and $TM$ given by the symplectic form $\o$: to this end it suffices to consider the fibers at the origin, that are
\begin{eqnarray*}
    T_o N_M &=&\gn^- \subset T_oM = \gn^- + \gn^+\, ,\\
    T^+_o N_M&=&\gn^+\, ,
\end{eqnarray*}
respectively, and to observe that
\begin{eqnarray*}
    \gn^+&\longrightarrow &(\gn^{-})^*\, ,\\
    X^+ &\longrightarrow &\o_o(X^+,\, \cdot\, )\, ,
\end{eqnarray*}
is an isomorphism.\par
Therefore, the flag manifold $F$ can be replaced by the open dense subset $N_F\subset F$: taking into account  the $N^-$--equivariant isomorphism $T^*N_F\simeq T^+M|_{N_M}$, it remains to construct a diffeomorphism $\Phi: T^+M|_{N_M} \longrightarrow M $ onto $M$. This will be given by the following exponential map:
\begin{eqnarray} 
\Phi:  T^+M|_{N_M}& \longrightarrow& M \, ,\nonumber\\
        n^-_*X^+ &\longmapsto& n^-(\exp X^+)o,\,\label{eqMappaDarbouxDmitri}
\end{eqnarray}
where   $n^- \in  N^-$, $X^+\in T_0^+M=\gn^+$,    $\exp : \gn^+ \to N^+$ is the  exponential map  of the Lie  algebra $\gn^+$  into the  the group $N^+$, and $n^-_*:T^+_oM\to T^+_{n^-o}M$ is the differential at $o$ of the action of $n^-$.\par
By  Corollary~\ref{corIntegralsOrbits}, the restriction $\Phi|_{{T^+_{n^-0}}M}$ of $\Phi$ to the fiber $T^+_{n^-0}M$ of the rank--$n$ bundle $ T^+M|_{N_M}$ at the point $n^-o$ of the $n$--dimensional manifold $N_M$ is a diffeomorphism onto its image, which is the $n$--dimensional maximal integral submanifold of $T^+M$ that passes through $n^-o$; since the set of all points of the form $n^-o$, that is the orbit $N^-o=N_M$, is transversal to the aforementioned integral manifolds, the latter make up a (smooth) bundle over $N_M$: but~\eqref{eqMappaDarbouxDmitri} is manifestly $N^-$--equivariant, so that the whole map $\Phi$, regarded as the $N^-$--equivariant extension of the diffeomorphism $\Phi|_{{T^+_{n^-0}}M}$, is a diffeomorphism itself.
It remains to show that the diffeomorphism~\eqref{eqMappaDarbouxDmitri}, regarded as a map
\begin{eqnarray} 
\Phi:  T^*N_F& \longrightarrow& M \, ,\nonumber
\end{eqnarray}
pulls back the symplectic form $\omega_Z$ given by~\eqref{eqSympFormAdZ} to the ``deformation'' $\omega_{T^*F}$ of the standard symplectic form $\omega_{\mathrm{std}}$ on $T^*N_F$ given by~\eqref{eqDeformedSymplecticForm}. 
To begin with, we compare $\omega_{T^*F}$ with $\Phi^*(\omega_Z)$ at the zero $0\in T_{\pi(o)} ^*N_F$ of the fiber of $T^*N_F$ at the origin $\pi(o)$ of $N_F$.\par 
Let us recall that
\begin{equation}
\omega_{\mathrm{std}}|_0\in\Lambda^2(T^*_0(T^*N_F))=\Lambda^2((T_{\pi(o)}N_F\oplus T_0(T_{\pi(o)}^*N_F))^*)=\Lambda^2(T^*_{\pi(o)}N_F\oplus T_{\pi(o)}N_F)
\end{equation}
is given by
\begin{equation}
    \omega_{\mathrm{std}}|_0=e_i\wedge\varepsilon^i\, ,
\end{equation}
where $\{e_i\}$ is a basis of $T_{\pi(o)}N_F$ and $\{\varepsilon^i\}$ is its dual basis in $T^*_{\pi(o)}N_F$.\par 
We shall need a root system $R$  of $\gg$, where $R=R^+\cup  R^-$ is its stplitting into positive and negative roots, and $\gg_\alpha$ denotes the eigenspace of each $\alpha\in R$.\par  
Then it is possible to choose as basis of $T_{\pi(o)}N_F\equiv\gn^-$ a system $\{E_{\alpha}\mid \alpha\in R^+\}$, where each $E_\alpha$ is a generator of $\gg_\alpha$ for all  $\alpha\in R$, and 
\begin{equation}
    (E_\alpha,E_{-\beta})=\delta_{\alpha,\beta}\, ,\quad\forall\alpha,\beta\in R^+\, .
\end{equation}
Since $\Phi$ maps $0\in T^*N_F$ to $o=Z\in M$, one needs to calculate
\begin{equation}
    \omega_Z(E_\alpha,E_\beta)
\end{equation}
for all $\alpha,\beta\in R$: by definition,
\begin{equation}
    \omega_Z(E_\alpha,E_\beta)=B(Z,[E_\alpha,E_\beta])\, ,
\end{equation}
and $[E_\alpha,E_\beta]=0$ unless $\alpha+\beta=0$, in which case
\begin{equation}
    \omega_Z(E_\alpha,E_{-\alpha})=B(Z,H_\alpha)=j\, ,
\end{equation}
where $j$ is the degree of $E_\alpha$, 
which matches with the definition~\eqref{eqDeformedSymplecticForm} of $\omega_{T^*F}$.\par
The general claim follows from invariance arguments.
\end{proof}

\section{Monge--Amp\`ere equations and their    generalization}\label{secMAEs}

In this section we review some basics fact concerning the contact geometry of Monge--Amp\`ere equations. More details can be found in \cite{MR2985508,MR2352610}.

\subsection{Monge--Amp\`ere equations via $n$--forms on jet spaces}
With any differential form $\Omega$ on the $(2n+1)$--dimensional jet space 
\begin{equation}
    J^1(n,1):=\{[f]_x^1\mid f:\R^n\to\R\}
\end{equation}
one can associate its so--called \textit{horizontalization}: such an operation corresponds to the projection of $\Lambda^\bullet(J^1)$ onto the quotient 
\begin{equation}
 \overline{\Lambda}^\bullet(J^1)\df \frac{\Lambda^\bullet(J^1)}{\mathcal{I}_{\mathcal{C}}}\, ,
\end{equation}
where  $\Lambda^\bullet(J^1)$ is the algebra of differential forms and ${\mathcal{I}_{\mathcal{C}}}$ is the differential ideal  generated by the contact form
\begin{equation}
    du-\sum_{i=1}^nu_idx^i\, .
\end{equation}
 We denote by $\overline{\Omega}$ the image of $\Omega\in \Lambda^\bullet(J^1)$ in $\overline{\Lambda}^\bullet(J^1)$.\par
It is not hard to see that
\begin{equation}
     \overline{\Lambda}^n(J^1)=C^\infty(J^1)\otimes_{C^\infty(\R^n)}\Lambda^n(\R^n)\, ,
\end{equation}
that is, for any $\Omega\in \Lambda^n(J^1)$ there exists a unique function $F_\Omega\in C^\infty(J^1)$, such that $\overline{\Omega}=F_\Omega\cdot \dd x^1\wedge\cdots\wedge\dd x^n$.
\begin{definition}\label{defMAE1}
    For any $n$--form $\Omega$ on $J^1$, we call $\E_\Omega :=\{F_\Omega=0\}$ the \emph{Monge--Amp\`ere equation (MAE)} associated with $\Omega$ (see Section~\ref{secMAL_contact} below). A function $f:\mathcal{U}\subseteq\R^n\to\R$, such that $[f]_x^2\in\E_\Omega$   $\forall x\in \mathcal{U}$, is a \emph{solution} of $\E_\Omega$, whereas each point $[f]_x^2\in\E_\Omega$ is called a \emph{formal solution} of $\E_\Omega$
\end{definition}

\subsection{Monge--Amp\`ere equations on  contact  and  symplectic manifolds}\label{secMAL_contact}

Given a {contact manifold} $(N,\mathcal{C})$, that is a $(2n+1)$--dimensional manifold equipped with a completely non--integrable $2n$--dimensional distribution $\mathcal{C}$, formula $\omega:=d\theta|_{\mathcal{C}}$ defines a conformal symplectic structure  on $\mathcal{C}$, where $\theta$ is any $1$--form, such that $\ker(\theta)=\mathcal{C}$. In turn, this conformal symplectic structure allows to introduce
 the \emph{Lagrangian Grassmannian} of $(\mathcal{C}_p, \omega_p)$, that is the set
\begin{eqnarray*}
    \mathcal{L}(\mathcal{C}_p)&\df&\{L_p\mid L_p\text{ is a Lagrangian planes of }\mathcal{C}_p\}\\
    &=&\{L_p\in\Gr(n,\mathcal{C}_p)\mid\omega_p|_{L_p}\equiv 0\}\, 
\end{eqnarray*}
of all Lagrangian planes at a given point   $p\in N$.
\begin{definition}\label{def.M1}
The \textit{prolongation} of a contact manifold $(N,\mathcal{C})$ is the fiber bundle $\pi: N^{(1)}\longrightarrow N$, where
$$
N^{(1)}:=\bigcup_{p\in N}\mathcal{L}(\mathcal{C}_p) \, ,
$$
and $\pi$ is the natural projection.
\end{definition}
Points $p^1$ of $N^{(1)}$ can be then understood as Lagrangian planes $L_{p^1}$ of $(\mathcal{C}_p,\omega_p)$, that is, there is a  natural correspondence
$$
p^1\in N^{(1)}\Longleftrightarrow L_{p^1}\in\mathcal{L}(\mathcal{C}_p)\, ,\quad p=\pi(p^1)\,.
$$
According to V.~Lychagin \cite{MR2352610}  and T.~Morimoto \cite{Morimoto1979}, we generalize Definition~\ref{defMAE1} as follows. Given a  contact manifold $(N,\mathcal{C})$, with  $\mathcal{C}=\ker(\theta)$, let  $\mathcal{I}_\CC$ denote the differential ideal generated by $\theta$; recall also that   $L_{{p}^1}\subset T_{\pi({p}^1)}N$ is the
Lagrangian plane associated with ${p}^1\in N^{(1)}$.
\begin{definition}\label{def.MA.Lych}
The hypersurface $\mathcal{E}_{\Omega}$ of $N^{(1)}$ given by
\begin{equation}\label{eqDefMAE_general}
\mathcal{E}_{\Omega}:=\{{p}^1\in N^{(1)}
\mid \Omega|_{L_{{p}^1}}\equiv 0\}\, ,
\end{equation}
where $\Omega\in\Lambda^{n}(N)$, is called a \emph{(general)  Monge--Amp\`ere equation (MAE)}.
\end{definition}
MAEs can be defined also on symplectic manifolds rather than on contact ones: it is enough to replace the above contact manifold $N$ with a symplectic manifold $M$:  details of such a construction are omitted it, see, e.g.,~\cite{Ferapontov2020,Russo2019}.
\begin{remark}
Since, in the definition~\eqref{eqDefMAE_general}
of  a general MAE, the restriction $\Omega|_{L_{{p}^1}}$ is identically zero   if the form  $\Omega\in\Lambda^{n}(N) $ belongs to the contact ideal $\mathcal{I}_\CC$, instead of the form $\Omega$ one can use the equivalence class \begin{equation}  [\Omega]\in\frac{\Lambda^{n}(N)}{\mathcal{I}_\CC}\, .
\end{equation}
Elements of the quotient
\begin{equation}  \frac{\Lambda^{\bullet}(N)}{\mathcal{I}_\CC}
\end{equation}
are called  \emph{effective} forms; accordingly, some speak of MAE associated with an \emph{effective} differential form.
\end{remark}
\subsection{MAEs in Darboux coordinates}\label{secMAEinCoordinates}
If $F$ is an $n$--dimensional smooth manifold, equipped with coordinates $\{x^1,\ldots,x^n\}$, then the standard symplectic form $\omega_{\mathrm{std}}$ on $T^*F$ reads
\begin{equation}
    \omega_{\mathrm{std}}=dx^i\wedge du_i\, ,
\end{equation}
where $u_i$ is the momentum conjugate with $x^i$; together with the coordinate $u$ that corresponds to the value of a function, the $x^i$'s and the moments $u_i$'s form a Darboux coordinate system of the contact manifold
\begin{equation}
    J^1(F,\R)=T^*F\times\R\, .
\end{equation}
In particular, the contact form $\theta$ reads
\begin{equation}
    \theta=du-u_idx^i\, .
\end{equation}
The second--order jet space $J^2(F,\R)$ is an affine bundle over $J^1(F,\R)$, its   fiber at  $p^1\in J^1(F,\R)$ being modeled by the symmetric power $S^2(T_p^*F)$, where $p=\pi(p^1)$: this allows to extend the Darboux coordinate system by adding the coordinates $u_{ij}$ that correspond to second--order derivatives. A point $p^2$ of the  fiber of $J^2(F,\R)$ over $p^1\in J^1(F,\R)$ has coordinates $u_{ij}$, if
\begin{equation}
    u_{ij}dx^idx^j\in S^2(T_p^*F)
\end{equation}
is  the   symmetric form that corresponds to $p^2$; the same point can be regarded as a Lagrangian plane of $\CC_{p^1}$ by means of the embedding
\begin{eqnarray}
    S^2(T_p^*F) &\longrightarrow &\LL(\CC_{p^1})\, ,\nonumber \\
    p^2=u_{ij}dx^idx^j&\longmapsto & L_{p^2}=\Span{D_{x^i}+u_{ij}\partial_{u_{j}}\mid i=1,\ldots,n}\, ,\label{eqEmbeddingS2LGr}
\end{eqnarray}
where $D_{x^i}=\partial_{x^i}+u_i\partial_u$ are the total derivatives.\par
A coordinate expression of the general MAE $\E_{\Omega}$ defined by an $n$--form $\Omega$ on $J^2(F,\R)$ can be then obtained by employing the extended Darboux coordinates introduced above; to this end we shall need the Pl\"ucker embedding
\begin{eqnarray}
    \LL(\CC_{p^1})&\longrightarrow&\mathbb{P}(\Lambda^n\CC_{p^1})\, ,\nonumber\\
    L=\Span{\ell_1\ldots,\ell_n}&\longmapsto &\vol(L):=[\ell_1\wedge\cdots\wedge\ell_n]\, ,\label{eqPlucker}
\end{eqnarray}
as well as a coordinate expression
\begin{equation}\label{eqCoordExprExOmega}
    \Omega=Adx^1\wedge\cdots\wedge dx^n+B_{i}^jdx^1\wedge\cdots\wedge \widehat{dx^i}\wedge\cdots\wedge dx^n\wedge du_j+\cdots+Cdu_1\wedge\cdots\wedge du_n
\end{equation}
of $\Omega$. By combining~\eqref{eqEmbeddingS2LGr} with~\eqref{eqPlucker}, we see that the point $p^1$ is mapped to the projective class of
\begin{equation}
    D_{x^1}\wedge\cdots\wedge D_{x^n}+u_{ij} D_{x^1}\wedge\cdots\wedge \widehat{D_{x^i}}\wedge\cdots\wedge D_{x^n}\wedge  \partial_{u_j}+\cdots+\det(u_{ij})\partial_{u_1}\wedge\cdots\wedge\partial_{u_n}\, .
\end{equation}
Then, by applying the formula~\eqref{eqDefMAE_general} that defines a MAE, we find out that
\begin{equation}
    \E_\Omega=\{F_\Omega(u_{ij})=0\}\, ,
\end{equation}
where $F_\Omega(u_{ij})$ is a linear combination of the minors of the matrix $(u_{ij})$, unambiguously defined by
\begin{equation}\label{eqFunctionMAE}
    \Omega|_{L_{p^2}}=F_\Omega(u_{ij})dx^1\wedge\cdots\wedge dx^n\, .
\end{equation}
\begin{example}
    If all coefficients in~\eqref{eqCoordExprExOmega} are zero, except $C$, then $ \E_\Omega$ is the   Monge--Amp\`ere equation given by $F_\Omega(u_{ij})=C\det(u_{ij})$, that is $\det(u_{ij})=0$.
\end{example}
We refer the reader to~\cite{Gutt2019,MR3760967} and references therein for more details.
\subsection{Fibers of MAE as hyperplane sections}\label{secFiberWise}
At each point $p$ of  the $(2n+1)$--dimensional contact manifold $(N,\CC)$, let us consider the Pl\"ucker embedding~\eqref{eqPlucker}. Then, the fiber
\begin{equation}
(\E_\Omega)_p:=\E_\Omega\cap N^{(1)}_p=\E_\Omega\cap \LL(\CC_p)
\end{equation}
of the MAE $\E_\Omega$ defined by~\eqref{eqDefMAE_general} turns out to be a \emph{hyperplane section} of $\LL(\CC_p)$. Indeed, the evaluation $\Omega_p$ of $\Omega$ at $p\in N$ is an element of the exterior algebra
\begin{equation} \Lambda^n(\CC_p^*)=(\Lambda^n(\CC_p))^*
\end{equation}
of the dual vector space $\CC_p^*$ and, as such, the linear equation
\begin{equation}
    \Omega_p=0
\end{equation}
defines a projective hyperplane in $\mathbb{P}(\Lambda^n\CC_p)$, whose intersection with $\LL(\CC_p)$ gives precisely $(\E_\Omega)_p$.\par
In view of the natural projection
\begin{eqnarray}
T_p^*N\longrightarrow\CC_p^*
\end{eqnarray}
the basis 
\begin{equation}
    (dx^1)_p\, ,\ldots\, , (dx^n)_p\, , (du)_p\, , (du^1)_p\, ,\ldots\, ,(du^n)_p\, ,
\end{equation}
of $T_p^*N$ gives rises to the basis
\begin{equation}
    (dx^1)_p\, ,\ldots\, , (dx^n)_p\, ,   (du^1)_p\, ,\ldots\, ,(du^n)_p\, , \label{eqBasePuntualeCiPiStar}
\end{equation}
of $\CC_p^*$. In the present paper, we will use~\eqref{eqBasePuntualeCiPiStar} as a standard basis of  $\CC_p^*$; moreover, since homogeneity  allows to restrict ourselves to the fiber of $N$ at the origin, the index $p$ will be omitted in the symbols above. 
\section{Examples of para--K\"ahler homogeneous manifolds and their contactifications}\label{secExamplesPK}


\subsection{Fundamental gradations of a (complex or real) semisimple Lie algebra $\gg$}
Let 
\begin{equation}
    \gg=\CSUB+\sum_{\alpha\in R}\gg_\alpha
\end{equation}
be a root space decomposition of a complex semisimple Lie algebra $\gg$ with respect to a Cartan subalgebra $\CSUB$. We fix a system of simple roots $\Pi = \{\alpha_1, \ldots , \alpha_\ell\} \subset R$, that is a basis of $\CSUB^*$, such that any root $\alpha \in R$ has integer coefficients with respect to $\Pi$ of the same sign (non--negative or non--positive).\par 
Any disjoint decomposition 
\begin{equation}
    \Pi=\Pi^0\cup\Pi^1
\end{equation}
of $\Pi$ defines a fundamental gradation of $\gg$ as follows. First, define the function $d : R \longrightarrow\Z$ by 
\begin{eqnarray*}
    d|_{\Pi^0} &:=& 0\, ,\\
     d|_{\Pi^1} &:=& 1\, ,\\
      d(\alpha) &:=& \sum k_id(\alpha_i)\, ,\quad \forall \alpha=\sum k_i\alpha_i\in R\, .
\end{eqnarray*} 
Then the fundamental gradation is given by
\begin{equation}
    \gg^0=\CSUB+\sum_{\alpha\in R_0} \gg_{\alpha}\, ,\quad R_0:=\{\alpha\in R\mid d(\alpha)=0\}\, ,
\end{equation}
and
\begin{equation}
    \gg^i=\sum_{\alpha,\,    d(\alpha)=i} \gg_{\alpha}\, , \quad \forall i\neq 0\, .
\end{equation}
Notice that any fundamental gradation of $\gg$ is conjugated to a unique gradation of such a form.\par 
Any real semisimple Lie algebra $\widehat{\gg}$ is a real form of a complex semisimple Lie algebra $\gg$, that is $\widehat{\gg}=\gg^{\sigma}$ is the fixed point set  of some antilinear involution $\sigma$ of $\gg$.  
We can always assume that $\sigma$ preserves a Cartan subalgebra $\CSUB$ of $\gg$ and induces an automorphism of the root system $R$. A root $a \in R$ is called \textit{compact} (or \textit{black}) if $\sigma(\alpha) =-\alpha$. It is always possible to choose a system of simple roots $\Pi = \{\alpha_1,\ldots,\alpha_\ell\}$ such that, for any non--compact root $\alpha_i \in \Pi$, the corresponding root $\sigma(\alpha_i)$ is a sum of one non--compact root $a_j \in \Pi$ and a linear combination of compact roots from $\Pi$. The roots $\alpha_i$ and $\alpha_j$ are called \textit{equivalent}.\par 
\begin{proposition}[Alekseevsky--Medori~\cite{Alekseevsky20078}] Let $\gg$ be a complex semisimple Lie algebra, $\sigma : \gg \longrightarrow\gg $ an antilinear involution, and $\gg^\sigma$ the corresponding real form. The gradation of $\gg$, associated with a decomposition $\Pi = \Pi^0 \cup \Pi^1$, defines a gradation $\gg^\sigma= \Sigma(\gg^i)^\sigma$ of $\gg^\sigma$ if and only if $\Pi^1$ consists of non--compact roots and any two equivalent roots are either both in $\Pi^0$ or both in $\Pi^1$. 
\end{proposition}
\subsection{Examples of fundamental gradations}
\subsubsection{Fundamental gradations of $\sll(V)$}
Let $V$ be a (complex or real) vector space and $V = V^1+\cdots+ V^k$ a decomposition of $V$ into a direct sum of subspaces. It defines a fundamental gradation 
\begin{equation}
    \sll(V)=\sum_{i=-k}^k\gg_i
\end{equation}
of the Lie algebra $\sll(V)$, where 
\begin{equation}\label{eqFundGradSLV}
    \gg^i = \{A \in \sll(V)\mid  A(V^j) \subset V^{i+j}\,\quad\forall   j = 1, \ldots , k\}\,  .
\end{equation}
\begin{proposition}
    Any fundamental gradation of $\sll(V)$ is of the form~\eqref{eqFundGradSLV}. 
\end{proposition}
\subsubsection{Fundamental gradations of $\gg_2$}\label{secFundGradG2}
The root system of the complex exceptional Lie algebra $\gg_2$ has the form 
\begin{equation}\label{eqRootSystemG2}
    R = \{\pm\epsilon_i,\pm(\epsilon_i-\epsilon_j)\mid i,j =1, 2, 3\}\, ,
\end{equation}
where the vectors $\epsilon_i$ satisfy 
\begin{eqnarray}
    \epsilon_1+\epsilon_2+\epsilon_3&=&0\, ,\label{eqRootSystemG2_1}\\
    \epsilon_i^2&=&\frac{2}{3}\, ,\label{eqRootSystemG2_2}\\
    (\epsilon_i,\epsilon_j)&=&-\frac{1}{3}\, ,\quad i\neq j\, .\label{eqRootSystemG2_3}
\end{eqnarray}
Consider the system of simple roots $\Pi = \{\alpha_1 := -\epsilon_2, \alpha_2 := \epsilon_2-\epsilon_3 \}$. The corresponding system of positive roots is 
\begin{equation}
    R^+ = \{\alpha_1, \alpha_2, \alpha_1+\alpha_2, 2\alpha_1 + \alpha_2, 3\alpha_1 + \alpha_2, 3\alpha_1 + 2\alpha_2\}\, ,
\end{equation}
and it is represented by the red arrows in the figure below:\par
\bigskip
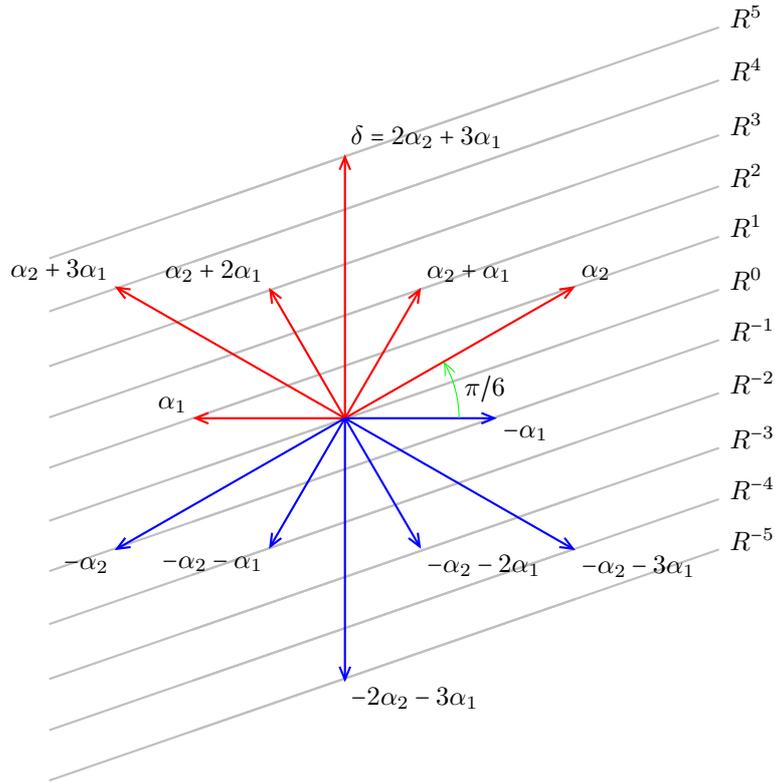
\begin{center}
      \begin{tikzpicture}[
    -{Straight Barb[bend,
       width=\the\dimexpr10\pgflinewidth\relax,
       length=\the\dimexpr12\pgflinewidth\relax]},
  ]
\label{fig}
  \draw [-,color={gray!50},thick] (3*60:2) -- (-30 + 1*60:3.5)--([turn]0:2);
  \draw [-,color={gray!50},thick] (-30 + 1*60:3.5)--(3*60:2) -- ([turn]0:2);

  \node[above right, inner sep=.2em] at (1*25:5.5) {$R^1$};

 \draw [-,color={gray!50},thick,yshift=19] (3*60:2) -- (-30 + 1*60:3.5)--([turn]0:2);
  \draw [-,color={gray!50},thick,yshift=19] (-30 + 1*60:3.5)--(3*60:2) -- ([turn]0:2);

    \node[above right, inner sep=.2em,yshift=19] at (1*25:5.5) {$R^2$};

  \draw [-,color={gray!50},thick,yshift=38.3] (3*60:2) -- (-30 + 1*60:3.5)--([turn]0:2);
  \draw [-,color={gray!50},thick,yshift=38.3] (-30 + 1*60:3.5)--(3*60:2) -- ([turn]0:2);

  \node[above right, inner sep=.2em,yshift=38.3] at (1*25:5.5) {$R^3$};

    \draw [-,color={gray!50},thick,yshift=59] (3*60:2) -- (-30 + 1*60:3.5)--([turn]0:2);
  \draw [-,color={gray!50},thick,yshift=59] (-30 + 1*60:3.5)--(3*60:2) -- ([turn]0:2);

    \node[above right, inner sep=.2em,yshift=59] at (1*25:5.5) {$R^4$};

  \draw [-,color={gray!50},thick,yshift=79] (3*60:2) -- (-30 + 1*60:3.5)--([turn]0:2);
  \draw [-,color={gray!50},thick,yshift=79] (-30 + 1*60:3.5)--(3*60:2) -- ([turn]0:2);

    \node[above right, inner sep=.2em,yshift=79] at (1*25:5.5) {$R^5$};


  \draw [-,color={gray!50},thick,yshift=-118] (3*60:2) -- (-30 + 1*60:3.5)--([turn]0:2);
  \draw [-,color={gray!50},thick,yshift=-118] (-30 + 1*60:3.5)--(3*60:2) -- ([turn]0:2);

  \node[above right, inner sep=.2em,yshift=-118] at (1*25:5.5) {$R^{-5}$};

 \draw [-,color={gray!50},thick,yshift=19-118] (3*60:2) -- (-30 + 1*60:3.5)--([turn]0:2);
  \draw [-,color={gray!50},thick,yshift=19-118] (-30 + 1*60:3.5)--(3*60:2) -- ([turn]0:2);

    \node[above right, inner sep=.2em,yshift=19-118] at (1*25:5.5) {$R^{-4}$};

  \draw [-,color={gray!50},thick,yshift=38.3-118] (3*60:2) -- (-30 + 1*60:3.5)--([turn]0:2);
  \draw [-,color={gray!50},thick,yshift=38.3-118] (-30 + 1*60:3.5)--(3*60:2) -- ([turn]0:2);

  \node[above right, inner sep=.2em,yshift=38.3-118] at (1*25:5.5) {$R^{-3}$};

    \draw [-,color={gray!50},thick,yshift=59-118] (3*60:2) -- (-30 + 1*60:3.5)--([turn]0:2);
  \draw [-,color={gray!50},thick,yshift=59-118] (-30 + 1*60:3.5)--(3*60:2) -- ([turn]0:2);

    \node[above right, inner sep=.2em,yshift=59-118] at (1*25:5.5) {$R^{-2}$};

  \draw [-,color={gray!50},thick,yshift=79-118] (3*60:2) -- (-30 + 1*60:3.5)--([turn]0:2);
  \draw [-,color={gray!50},thick,yshift=79-118] (-30 + 1*60:3.5)--(3*60:2) -- ([turn]0:2);

    \node[above right, inner sep=.2em,yshift=79-118] at (1*25:5.5) {$R^{-1}$};

 \draw [-,color={gray!50},thick,yshift=-20] (3*60:2) -- (-30 + 1*60:3.5)--([turn]0:2);
  \draw [-,color={gray!50},thick,yshift=-20] (-30 + 1*60:3.5)--(3*60:2) -- ([turn]0:2);

    \node[above right, inner sep=.2em,yshift=-20] at (1*25:5.5) {$R^{0}$};
    
    \foreach \i in {1, 2, 3} {
      \draw[thick, red] (0, 0) -- (\i*60:2);
      \draw[thick, red] (0, 0) -- (-30 + \i*60:3.5);
    }
    \foreach \i in {4, 5, 6} {
      \draw[thick, blue] (0, 0) -- (\i*60:2);
      \draw[thick, blue] (0, 0) -- (-30 + \i*60:3.5);
    }
    \draw[thin, green] (1.5, 0) arc[radius=1.5, start angle=0, end angle=1*30];
    \node[above left, inner sep=.2em] at (6*30:2) {$\alpha_1$};
    \node[below right, inner sep=.2em] at (0*30:2) {$-\alpha_1$};
    \node[above right, inner sep=.2em] at (1*30:3.5) {$\alpha_2$};
    \node[below left, inner sep=.2em] at (7*30:3.5) {$-\alpha_2$};
    \node[above right, inner sep=.2em] at (2*30:2) {$\alpha_2+\alpha_1$};
    \node[below left, inner sep=.2em] at (8*30:2) {$-\alpha_2-\alpha_1$};
    \node[above right, inner sep=.2em] at (3*30:3.5) {$\delta=2\alpha_2+3\alpha_1$};
    \node[below right, inner sep=.2em] at (9*30:3.5) {$-2\alpha_2-3\alpha_1$};
    \node[above left, inner sep=.2em] at (4*30:2) {$\alpha_2+2\alpha_1$};
    \node[below right, inner sep=.2em] at (10*30:2) {$-\alpha_2-2\alpha_1$};
    \node[above left, inner sep=.2em] at (5*30:3.5) {$\alpha_2+3\alpha_1$};
    \node[below right, inner sep=.2em] at (11*30:3.5) {$-\alpha_2-3\alpha_1$};
    \node[right] at (15:1.5) {$\pi/6$};
  \end{tikzpicture}
\end{center}
\bigskip
There are three fundamental gradations for the complex Lie algebra $\gg_2$ . For any of such gradations, we give below the subset $\Pi^1 \subset\Pi$ and the level sets $R^i := \{\alpha \in R\mid d(\alpha) = i\}$ of the grading function $d: R\longrightarrow\Z$.
\begin{itemize}
    \item[1)] $\Pi^1=\Pi$:
\begin{eqnarray*}
    R^0&=&\emptyset\, ,\\
    R^1&=&\{\alpha_1,\alpha_2\}\, ,\\
    R^2&=&\{\alpha_1+\alpha_2\}\, ,\\
    R^3&=&\{2\alpha_1+\alpha_2\}\, ,\\
    R^4&=&\{3\alpha_1+\alpha_2\}\, ,\\
    R^5&=&\{3\alpha_1+2\alpha_2\}\, .
\end{eqnarray*} 
All the level sets $R^i$ , with $i=-5,\ldots,+5$ are represented by gray parallel lines in the picture above.
    \item[2)] $\Pi^1 = \{\alpha_1\}$: 
    \begin{eqnarray*}
    R^0&=&\{\alpha_2\}\, ,\\
    R^1&=&\{\alpha_1,\alpha_1+\alpha_2\}\, ,\\
    R^2&=&\{2\alpha_1+\alpha_2\}\, ,\\
    R^3&=&\{3\alpha_1+\alpha_2,3\alpha_1+2\alpha_2\}\, .
\end{eqnarray*} 
    \item[3)]$\Pi^1 = \{\alpha_2\}$: 
    \begin{eqnarray*}
    R^0&=&\{\alpha_1\}\, ,\\
    R^1&=&\{\alpha_2,\alpha_1+\alpha_2,2\alpha_1+\alpha_2,3\alpha_1+\alpha_2\}\, ,\\
    R^2&=&\{3\alpha_1+2\alpha_2\}\, .
\end{eqnarray*} 
\end{itemize}
There are just two real forms of the complex Lie algebra $\gg_2$: the compact form, which has no non--trivial gradation, and the normal form $\gg_2^\sigma$, which has a diagonalizable Cartan subalgebra and no compact roots. The above--listed gradations of the complex Lie algebra $\gg_2$ define three gradations of the real Lie algebra $\gg_2^\sigma$.

\section{Classification of $\mathsf{G}_2$--invariant MAEs on the contactification of $\mathsf{G}_2/\GL_2(\R)$}\label{secG2_case}

\subsection{Invariant effective $n$--forms}


A key step towards a classification of invariant Monge--Amp\`ere equations on the contactification $\Mcont=G/L$ of a $G$--homogeneous para--K\"ahler manifold $M=G/H$ (see Section~\ref{secContactification}) is the description of  effective   $G$-invariant   $n$--forms  $\Omega \in \Lambda^n(M)^G$  on the $2n$--dimensional symplectic manifold   $M$: in turn, this reduces  to the  description of the  space $\L^n(\gn^*)^H$, which is a rather  simple  problem, since    the  $n$--form  must be invariant with respect to  the  diagonal  operator $\ad_Z$.\par

Then, in order to describe $G$--invariant Monge--Amp\`ere equations on $\Mcont$, it suffices to describe the  $\ad_{\gl}$--invariant  $n$--forms on  $\gn$; moreover,  the problem itself can be simplified, if one   assumes that    these forms   are    eigen--forms   for  $\ad_Z$.
%
%
\subsubsection{The space of $H$--invariant effective $n$--forms on $M$}
In view of the reductive  decomposition
\begin{equation}
     \gg  =  \gh + \gn = (\gl + \bR Z)+ \gn^+ + \gn^-\, ,
\end{equation}
 the problem  reduces    to  the   description  of the  space $ \L^n(\gn^*)^{H}$
 of  $\Ad_H$--invariant (or   equivalently, if $H$ is connected,   $\ad_{\gh}$--invariant )  exterior  $n$--forms
 in   the    space $\gn = \gn^+ + \gn^-$.\par
   Denote  by
\begin{equation}
    \gn^{*} = \sum_a(\gn^{\pm})^*_{ a}
\end{equation}
      the  $\ad_Z$--eigenspace decomposition  of  the   space $(\gn^{\pm})^*$, where $a$ indicates  the   eigenvalue  of  the eigenspace $(\gn^{\pm})^*_{ a}$.
      Then a    decomposable   $n$--form
\begin{equation}
     \Omega = \xi_{a_1}\wedge \xi_{a_2}\wedge\cdots \wedge \xi_{a_n}\, , 
\end{equation}
   where  $\xi_a \in \gn_\alpha^*$,  is $\ad_Z$--invariant  if and only if
\begin{equation}
    \sum a_i =0\, . 
\end{equation}
The description of  all invariant $n$--forms  on  $\gn$ is based on these simple remarks.
\subsubsection{The space of $L$--invariant effective $n$--forms on $\Mcont$}
    Since $\gn^{\pm}$  are   $\gl$-submoduli,   $\L^n(\gn)^L$ can be decomposed as
\begin{equation}
    \L^n(\gn)^L = \sum_{p+q=n} (\L^L)^{p,q} \, ,
\end{equation}
where
\begin{equation}
(\L^L)^{p,q}:= \L^p(\gn^+)^{L\,\ast}\wedge \L^q (\gn^-)^{L\,\ast}\, .
\end{equation}
It is then  sufficient to    describe $L$--invariant forms   in the  spaces
\begin{equation}
     \L^p(\gn^+)^{L\,\ast}\wedge \L^q (\gn^-)^{L\,\ast}\, .
\end{equation}
Note that  $\L^n(\gn^{\pm})^* \subset \L^n(\gn^*)^H $ if  $\Ad_H$ is unimodular.

\subsection{Properties of the roots of the exceptional Lie algebra  $\gg_2$}

We will need the formulas
\begin{eqnarray*}
    (E_\a, E_\a) &=& \frac{2}{(\a,\a)} \, , \\ \vspace{3mm} [H_{\a_i}, E_{\b}] &=&   \frac{2       (a_i,\b)}{(\a_i, \a_i)}E_{\b}\, ,
\end{eqnarray*}
where  $\a_i \in \Pi$ are  simple roots, $E_\alpha$ is a generator of $\gg_\alpha$,  and   $H_{\a} = [E_{\a}, E_{-\a}]$ for  any  root $\a \in R$: these are  rather standard; see, for instance,   Gorbatsevich, Vinberg and Onishchik's book~\cite{Onishchik1994-rv}.\par
Let  $\gg_2 $ be the   non--compact   real Lie  algebra   of type   $\mathsf{G}_2$ with the Cartan subalgebra  $\ga$; following the   notation  of the aforementioned book, we take the root system $R$ given by~\eqref{eqRootSystemG2}, where 
$\e_1$, $\e_2$ and  $\e_3$ satisfy \eqref{eqRootSystemG2_1}--\eqref{eqRootSystemG2_2}--\eqref{eqRootSystemG2_3}, and the symple roots are $\Pi = \{  \a_1 = -\e_2,\, \a_2 =\e_2 - \e_3\}$.
 with
\begin{equation}
     \a_1^2 =2/3\, ,\quad \a_2^2 =2\, ,\quad(\a_1, \a_2) =1\, ,
\end{equation}
see Section~\ref{secFundGradG2}.\par 
The maximal root  is
\begin{equation}
    \delta := 3 \a_1 + 2 \a_2 = \e_1 - \e_3\, ,
\end{equation}
in particular,
\begin{eqnarray}
    <\a|\delta>:= \frac{2<\a, \delta>}{\delta,
\delta}  &=& (\delta, \alpha)\, \quad\forall \a \in R\, ,\\
< \a_1| \delta> &=&0\, ,\\
< \a_2 + k \a_1|\delta>&=& (\a_2, \delta) =1\, \quad \forall k\, ,
\end{eqnarray}
%
see the picture at page~\pageref{fig}.\par 
The  operator $\ad_{H_{\delta}}$   acts  on   the  vector $E_{\a_2 + k \a_1}$ as  the identity $\Id$ and on the dual  form  $E^*_{\a_2 + k \a_1}$
as $- \Id$.\par 

\subsection{Invariant five--forms on the ten--dimensional manifold  $M  = \mathsf{G}_2/\GL_2(\R)$}
This case corresponds to the fundamental gradation 3) introduced earlier in Section~\ref{secFundGradG2}: the subalgebra $\gh$ of maximal rank  associated  with  the    root $\a_1$   is given   by
\begin{equation}
    \gh = \gl_2^{\a_1}= \bR H_{\delta} + \gsl_2(\R) = \R H_\delta  + \Span{H_{\a_1}, E_{\pm \a_1}}\, ,
\end{equation} 
and the corresponding subgroup $H$  of the  non compact real   exceptional group $\mathsf{G}_2$ will be given by
 $H =\GL_2(\R)$.
 Since  $(\a_1, \delta)=0$, the subgroup $H$ turns out ot be  the  centralizer    of the    element  $ Z=H_{\delta}$, 
 which acts   on   root  vectors   as
\begin{eqnarray}
    \ad_Z  E_{\a_2 +  k \a_1} &=& <\a_2 + k \a_1| \delta> E_{\a_2 + k \a_1} = E_{\a_2 + k \a_1}\, ,\\
 \ad_Z E_{\delta} &=& 2 E_{\delta}  \, ,    
\end{eqnarray}
 and   similarly  for  root vectors  corresponding to negative  roots. 
Therefore, the     reductive decomposition~\eqref{eqGenGaussAlg} reads   
\begin{equation*}
    \gg_2 = \gh + \gm = \gh + \gn_+ + \gn_-   \, ,
\end{equation*}
whereas the decomposition~\eqref{eqDecGiPiuMeno} reduces to
\begin{equation}
     \gn_{\pm} = \gn_{\pm 1} + \gn_{\pm 2}\, .   \label{eqEnnePiuMenoEnneUnoEnneDue}
\end{equation}
 The  tangent space  $\gm = T_o M$ of the   manifold   $M = \mathsf{G}_2/H =\mathsf{G}_2/ \GL_2(\bR) $    have the  basis
$\{E_{\pm \gamma_i}, E_{\pm\delta}\mid  i=0,1,2,3\}$,    where
\begin{equation}\label{eqGammaBase10DManG2}
    \gamma_i :=\a_2 +  i\a_1\, ,\quad  i=0,1,2,3\, ,\quad \delta= 2 \a_2 +3 \a_1  \, ,
\end{equation}
 are the positive roots    of  $\gg_2$, 
whereas  $\{H_{\delta}, H_{\a_1}, E_{\pm \a_1}\}$ is a basis of the stability  subalgebra  $\gh $, and  $\{ H_{\a_1},  E_{\pm \a_1}\}$ is a basis of the corresponding  derived subalgebra $\gh' = \gsl_2(\bR) $.\par 
In view of~\eqref{eqGammaBase10DManG2}, 
the action  of the operator  $\ad_{E_{\pm\a_1}}$   on  the root  vectors from $\gm^+$ is given by
\begin{eqnarray}
  \ad_{E_{\pm \a_1}} E_{\pm \delta}&=&0\, ,\label{eqActionEdelta}\\ \ad_{E_{\pm \a_1}} E _{\gamma_i}&=&
 N_{\pm\a_1,\gamma_i}E_{\gamma_{i\pm 1}}\, ,          \label{eqActionEgamma}
\end{eqnarray}
where the only nonzero $ N_{\a_1,\gamma_i}$'s are:
\begin{eqnarray*}
    N_{\a_1,\gamma_0}=N_{\a_1,\alpha_2}&=&1\, ,\\  N_{\a_1,\gamma_1}&=&2\, ,\\   N_{\a_1,\gamma_2}&=&3\, ,\\   
     N_{\a_1,-\gamma_1}&=&-3\, ,\\   N_{\a_1,-\gamma_2}&=&-2\, ,\\     N_{\a_1,-\gamma_3}&=&-1\, ,
\end{eqnarray*}
whereas the only nonzero  $ N_{-\a_1,\gamma_i}$'s are:
\begin{eqnarray*}
    N_{-\a_1,\gamma_1}&=&3\, ,\\  N_{-\a_1,\gamma_2}&=&2\, ,\\   N_{-\a_1,\gamma_3}&=&1\, ,\\    
    N_{-\a_1,-\gamma_0}=N_{-\a_1,-\alpha_2}&=&-1\, ,\\ N_{-\a_1,-\gamma_1}&=&-2\, ,\\   N_{-\a_1,-\gamma_2}&=&-3\, ,
\end{eqnarray*}
see, eg.,~\cite[Section~33.5]{humphreys2012linear}.
In the next propositions we will identify $\gm$ with its dual $\gm^*$.
%
\begin{proposition}\label{propOmega1}
    The space $\Lambda^1(\gm^*)^{\gh'}$ of $\gh'$--invariant one--forms on $\gm$ is generated by $E_{\pm\delta}^*$.
\end{proposition}
\begin{proof}
If
    \begin{equation}
        \omega^1=\mu_+E_{\delta}+\sum_{i=0}^3\lambda^i_+E_{\gamma_i}+\sum_{i=0}^3\lambda^i_-E_{-\gamma_i}+\mu_-E_{-\delta}\, ,
    \end{equation}
    then \eqref{eqActionEdelta}--\eqref{eqActionEgamma} show that
    \begin{equation}
        \ad_{E_{\alpha_1}}(\omega^1)= \sum_{i=0}^2\lambda^i_+N_{\alpha_1,\gamma_i}E_{\gamma_i}+\sum_{i=1}^3\lambda^i_-N_{\alpha_1,-\gamma_i}E_{-\gamma_i} 
    \end{equation}
    vanishes if and only if $\lambda_+^0=\lambda_+^1=\lambda_+^2=\lambda_-^1=\lambda_-^2=\lambda_-^3=0$, i.e., $\omega^1=\mu_+E_{\delta}+\lambda^3_+E_{\gamma_3}+\lambda^0_-E_{-\gamma_0}+\mu_-E_{-\delta}$. We can now apply $\ad_{E_{-\alpha_1}}$ and find out that 
    \begin{equation}
        \ad_{E_{-\alpha_1}}(\omega^1)=   \lambda^3_+N_{-\alpha_1,\gamma_3}E_{\gamma_3}+\lambda^0_-N_{-\alpha_1,-\gamma_0}E_{-\gamma_0}   
    \end{equation}
    vanishes if and only if $\lambda^3_+=\lambda^0_-=0$, whence  $\omega^1=\mu_+E_{\delta} +\mu_-E_{-\delta}$.
\end{proof}
\begin{proposition}\label{propOmega2}
    The space $\Lambda^2(\gm^*)^{\gh'}$ of $\gh'$--invariant two--forms on $\gm$ is generated by $E_{\delta}\wedge E_{-\delta}$, together with 
    \begin{eqnarray}
        \omega^2_\pm&:=&E_{\pm\gamma_1}\wedge E_{\pm\gamma_2}-3E_{\pm\gamma_0}\wedge E_{\pm\gamma_3}\, ,\label{eq_Inv_TWO_FORM}\\
        \omega^2&:=&  3 E_{\gamma_0}\wedge E_{-\gamma_0}+E_{\gamma_1}\wedge E_{-\gamma_1}+E_{\gamma_2}\wedge E_{-\gamma_2}+3E_{\gamma_3}\wedge E_{-\gamma_3}\,  .
    \end{eqnarray}
\end{proposition}

\begin{proof}
    From the splitting~\eqref{eqEnnePiuMenoEnneUnoEnneDue} it follows the splitting
    \begin{equation}
        \Lambda^2(\gm)=\Lambda^2(\gm_{+1})\oplus (\gm_{+1}\otimes \gm_{-1})\oplus  \Lambda^2(\gm_{-1}) \oplus(\gm_{+1}\otimes \gm_{+2})\oplus (\gm_{-1}\otimes \gm_{+2})\oplus (\gm_{+1}\otimes \gm_{-2})\oplus (\gm_{-1}\otimes \gm_{-2})\oplus (\gm_{+2}\otimes \gm_{-2}) 
    \end{equation}
    of  the 45--dimensional space $\Lambda^2(\gm)$. Since $\ad_{E_{\pm\alpha_1}}(\gm_{\pm 1})\subseteq \gm_{\pm 1}$ and $\ad_{E_{\pm\alpha_1}}(\gm_{\pm 2})=0$, all the constituents of the above decomposition are $\gh'$--invariant, so that we can analyze each of one separately. \par
We begin by observing that
\begin{equation}
        \ad_{E_{\pm\alpha_1}}(\omega^1\wedge E_{\pm\delta})= \ad_{E_{\pm\alpha_1}}(\omega^1)\wedge E_{\pm\delta}=0
    \end{equation}
    if and only if $\ad_{E_{\pm\alpha_1}}(\omega^1)=0$, i.e., $\omega^1$ is $\gh'$--invariant: by Proposition~\ref{propOmega1}, $\omega^1$ must then be a linear combination of $E_{\pm\delta}$, so that  the only invariant two--form in the last five constituents is $E_{ \delta}\wedge E_{-\delta}$.\par 
Passing now to $\Lambda^2(\gm_{+1})^{\gh'}$, it is convenient to introduce the basis
\begin{equation}
    \omega_{ij}:=E_{\gamma_i}\wedge E_{\gamma_j}\label{eqOmegaIJ}
\end{equation}
of the six--dimensional space $\Lambda^2(\gm_{+1})$: if
\begin{equation}\label{eqOmegaDUE}
    \omega^2=\lambda^{01}\omega_{01}+\lambda^{02}\omega_{02}+\lambda^{03}\omega_{03}+\lambda^{12}\omega_{12}+\lambda^{13}\omega_{13}+\lambda^{23}\omega_{23}\, ,
\end{equation}
then
\begin{eqnarray*}
    \ad_{E_{ \alpha_1}}(\omega^2)&=&\lambda^{01}\ad_{E_{ \alpha_1}}(\omega_{01})+\lambda^{02}\ad_{E_{ \alpha_1}}(\omega_{02})+\lambda^{03}\ad_{E_{ \alpha_1}}(\omega_{03})+\lambda^{12}\ad_{E_{ \alpha_1}}(\omega_{12})+\lambda^{13}\ad_{E_{ \alpha_1}}(\omega_{13})+\lambda^{23}\ad_{E_{ \alpha_1}}(\omega_{23})\\
    &=& \lambda^{01}N_{\alpha_1,\gamma_1} \omega_{02}+\lambda^{02} (N_{\alpha_1,\gamma_0}\omega_{12}+N_{\alpha_1,\gamma_2}\omega_{03})+\lambda^{03}N_{\alpha_1,\gamma_0}\omega_{13}+\lambda^{12} N_{\alpha_1,\gamma_2}\omega_{13}+\lambda^{13} N_{\alpha_1,\gamma_1}\omega_{23}
\end{eqnarray*}
vanishes if and only if
\begin{equation}
    \lambda^{03}N_{\alpha_1,\gamma_0} +\lambda^{12} N_{\alpha_1,\gamma_2}=\lambda^{03}  +3\lambda^{12}=0\, ,
\end{equation}
i.e., $\omega^2$ is proportional to $\omega_{12}-3\omega_{03}$: this two--form turn out to be $\ad_{E_{- \alpha_1}}$--invariant as well, since
\begin{equation*}
    \ad_{E_{- \alpha_1}}(\omega_{12}-3\omega_{03})=N_{-\alpha_1,\gamma_1}\omega_{02}-3N_{-\alpha_1,\gamma_3}\omega_{02}=(3-3\cdot 1)\omega_{02}=0\, .
\end{equation*}

The case of $\Lambda^2(\gm_{-1})^{\gh'}$ is formally analogous: instead of~\eqref{eqOmegaIJ} we shall have
\begin{equation}
    \omega_{ij}:=E_{-\gamma_i}\wedge E_{-\gamma_j}\,, 
\end{equation}
so that, having defined $\omega^2$ as in~\eqref{eqOmegaDUE}, it turns out that
\begin{eqnarray}
    \ad_{E_{ -\alpha_1}}(\omega^2) &=& \lambda^{01}N_{-\alpha_1,-\gamma_1} \omega_{02}+\lambda^{02} (N_{-\alpha_1,-\gamma_0}\omega_{12}+N_{-\alpha_1,-\gamma_2}\omega_{03})+\lambda^{03}N_{-\alpha_1,-\gamma_0}\omega_{13}\\
    &+&\nonumber\lambda^{12} N_{-\alpha_1,-\gamma_2}\omega_{13}+\lambda^{13} N_{-\alpha_1,-\gamma_1}\omega_{23}
\end{eqnarray}
vanishes if and only if
\begin{equation}
    \lambda^{03}N_{-\alpha_1,-\gamma_0} +\lambda^{12} N_{-\alpha_1,-\gamma_2}=-\lambda^{03}  -3\lambda^{12}=0\, ,
\end{equation}
i.e., $\omega^2$ is proportional to $\omega_{12}-3\omega_{03}$: this two--form turns out to be $\ad_{E_{ \alpha_1}}$--invariant as well, since
\begin{equation*}
    \ad_{E_{ \alpha_1}}(\omega_{12}-3\omega_{03})=N_{\alpha_1,-\gamma_1}\omega_{02}-3N_{\alpha_1,-\gamma_3}\omega_{02}=(-3-3\cdot(- 1))\omega_{02}=0\, .
\end{equation*}
By a further abuse of notation, we denote now by
\begin{equation}
    \omega_{ij}:=E_{\gamma_i}\wedge E_{-\gamma_j}
\end{equation}
the basis elements of the sixteen--dimensional space $\gm_{+1}\otimes \gm_{-1}$: if
\begin{equation}
    \omega^2=\sum_{i,j=0}^3\lambda^{ij}\omega_{ij}\, ,
\end{equation}
then
\begin{eqnarray*}
    \ad_{E_{ \alpha_1}}(\omega^2)&=&\sum_{i,j=0}^3\lambda^{ij}\ad_{E_{ \alpha_1}}(\omega_{ij})\\
    &=&\sum_{i,j=0}^3\lambda^{ij} (N_{\alpha_1,\gamma_i}\omega_{i+1,j}+N_{\alpha_1,-\gamma_j}\omega_{i,j-1})\\
    &=&\sum_{i,j=0}^3\lambda^{ij} N_{\alpha_1,\gamma_i}\omega_{i+1,j}+\sum_{i,j=0}^3N_{\alpha_1,-\gamma_j}\omega_{i,j-1}\, ,
\end{eqnarray*}
where the $\omega_{ij}$'s with an index beyond the range $\{0,1,2,3\}$ must be considered zero. Therefore
\begin{eqnarray*}
    \ad_{E_{ \alpha_1}}(\omega^2)&=& \sum_{j=0}^3\sum_{i=0}^2\lambda^{ij} N_{\alpha_1,\gamma_i}\omega_{i+1,j}+\sum_{i=0}^3\sum_{j=1}^3\lambda^{ij}N_{\alpha_1,-\gamma_j}\omega_{i,j-1} \\
    &=& \sum_{j=0}^3 (\lambda^{0j}\omega_{1j}+2\lambda^{1j}\omega_{2j}+3\lambda^{2j}\omega_{3j})  +\sum_{i=0}^3  ( -3 \lambda^{i1}\omega_{i0}    -2\lambda^{i2}\omega_{i1}    -\lambda^{i3}\omega_{i2}   )\\
    &=&\lambda^{00}\omega_{10}+2\lambda^{10}\omega_{20}+3\lambda^{20}\omega_{30}  +  \lambda^{01}\omega_{11}+2\lambda^{11}\omega_{21}+3\lambda^{21}\omega_{31} \\
    &&+\lambda^{02}\omega_{12}+2\lambda^{12}\omega_{22}+3\lambda^{22}\omega_{32} + \lambda^{03}\omega_{13}+2\lambda^{13}\omega_{23}+3\lambda^{23}\omega_{33} \\
    &&  -3 \lambda^{01}\omega_{00}    -2\lambda^{02}\omega_{01}    -\lambda^{03}\omega_{02}    -3 \lambda^{11}\omega_{10}    -2\lambda^{12}\omega_{11}    -\lambda^{13}\omega_{12} \\
    &&  -3 \lambda^{21}\omega_{20}    -2\lambda^{22}\omega_{21}    -\lambda^{23}\omega_{22}   -3 \lambda^{31}\omega_{30}    -2\lambda^{32}\omega_{31}    -\lambda^{33}\omega_{32} \\
    &=&   (\lambda^{00}-3\lambda^{11})\omega_{10}+(2\lambda^{10}-3\lambda^{21})\omega_{20}+3(\lambda^{20}-\lambda^{31})\omega_{30}\\
    &&+  (\lambda^{01}-2\lambda^{12})\omega_{11}+2(\lambda^{11}-\lambda^{22})\omega_{21}+(3\lambda^{21}-2\lambda^{32})\omega_{31}\\
    &&+ (\lambda^{02}-\lambda^{13} )\omega_{12}+(2\lambda^{12}-\lambda^{23})\omega_{22}+(3\lambda^{22}-\lambda^{33})\omega_{32} \\
    &&+ \lambda^{03}\omega_{13}+2\lambda^{13}\omega_{23}+3\lambda^{23}\omega_{33}-3 \lambda^{01}\omega_{00}    -2\lambda^{02}\omega_{01}    -\lambda^{03}\omega_{02}    \, .
\end{eqnarray*}
Analogously,
\begin{eqnarray*}
    \ad_{E_{ -\alpha_1}}(\omega^2)&=&  3\l^{10}\o_{00} +2\l^{20}\o_{10}+\l^{30}\o_{20} +(3\l^{11}-\l^{00})\o_{01}+(2\l^{21}-\l^{10})\o_{11}+(\l^{31}-\l^{20})\o_{21}\\
    &&+(3\l^{12}-2\l^{01})\o_{02}+2(\l^{22}-\l^{11})\o_{12}+(\l^{32}-2\l^{21})\o_{22}+3(\l^{13}-\l^{02})\o_{03}\\
    &&+(2\l^{23}-3\l^{12})\o_{13}+(\l^{33}-3\l^{22})\o_{23} -\l^{30}\o_{31}-\l^{31}\o_{32}-\l^{32}\o_{33}\, .
\end{eqnarray*}
It is not hard to see that the two equations $\ad_{E_{\pm \alpha_1}}(\omega^2)=0$ are satisfied if and only if $\lambda^{ij}=0$ for all $i\neq j$ and the following three conditions hold:
\begin{eqnarray*}
   \lambda^{11}&=&\lambda^{22}\, ,\\
   \lambda^{00}&=&3\lambda^{11}\, ,\\
   \lambda^{33}&=&3\lambda^{22}\, ,
\end{eqnarray*}
thus finishing the proof.
\end{proof}

\begin{proposition}\label{propFOURFORMS}
The space $\Lambda^4(\gm)^{\gh'}$ contain the following  (linearly independent) three four--forms:
    \begin{eqnarray}
    \o^4_{\pm} &:=& E_{\pm\gamma_0} \wedge E_{\pm\gamma_1}\wedge E_{\pm\gamma_2}\wedge E_{\pm \gamma_3 } \, ,\\
     \o^4 &:=&  E_{\gamma_0} \wedge E_{\gamma_1}\wedge E_{-\gamma_0}\wedge E_{-\gamma_1}+E_{\gamma_0} \wedge E_{\gamma_2}\wedge E_{-\gamma_0}\wedge E_{-\gamma_2}+\label{eq_MISSING_4form}\\
     &&+E_{\gamma_0} \wedge E_{\gamma_3}\wedge E_{-\gamma_1}\wedge E_{-\gamma_2}  +E_{\gamma_1} \wedge E_{\gamma_2}\wedge E_{-\gamma_0}\wedge E_{-\gamma_3} +\nonumber\\
     &&+E_{\gamma_1} \wedge E_{\gamma_3}\wedge E_{-\gamma_1}\wedge E_{-\gamma_3} +E_{\gamma_2} \wedge E_{\gamma_3}\wedge E_{-\gamma_2}\wedge E_{-\gamma_3}\, .\nonumber
\end{eqnarray}
\end{proposition}
\begin{proof}
    The six--dimensional $\sll_2$--representations $\Lambda^2(\gm_{\pm 1})$ can be decomposed into sums of irreducible $\sll_2$--representations, that are
    \begin{equation}
        \Lambda^2(\gm_{\pm 1})=\R\omega^2_\pm + V_4^{\pm}\, ,
    \end{equation}
    where $\omega^2_\pm$ are given by~\eqref{eq_Inv_TWO_FORM}, and $V_4^{\pm}$ are two copies of the irreducible $\sll_2$--representation of highest weight 4, that are dual to each other. The decomposition of the tensor product $\Lambda^2(\gm_{+ 1})\otimes\Lambda^2(\gm_{- 1})$ into  a sum of irreducible $\sll_2$--representations contains exactly two one--dimensional constituents: 
    \begin{equation}
        \Lambda^2(\gm_{+ 1})\otimes\Lambda^2(\gm_{- 1}) =\R\omega^2_+\wedge \omega^2_- + V_0+\textrm{irreducible modules of }\dim>1\, ,
    \end{equation}
    having employed the natural embedding
    \begin{equation*}
        \Lambda^2(\gm_{+ 1})\otimes\Lambda^2(\gm_{- 1})\subset \Lambda^4(\gm_{+ 1}\oplus\gm_{- 1})\, .
    \end{equation*}
    The one--dimensional representation $V_0$ is the unique one--dimensional constituent of $V_4^{+}\otimes V_4^{-}\simeq V_4 \otimes V_4^*$, which is is generated by 
    \begin{equation}
       \omega^4= \sum_{i=0}^4v_i\wedge v_i^*\, ,
    \end{equation}
    where $\{v_i\}$ and $\{v_i^*\}$ are bases of $V_4$ and $V_4^*$ dual to each other: easy computations lead to the  expression~\eqref{eq_MISSING_4form} of $\omega^4$ in our basis.\par
    The forms $\omega^4_{\pm}$ are the obvious generators of the one--dimensional  modules $\Lambda^4(m_{\pm 1})$.
\end{proof}

\begin{corollary}
    The algebra $\Lambda^\bullet(\gm)^{\gh'}$ of $\gh'$--invariant forms in generated, up to degree 5, by the eight elements
    \begin{equation}\label{eqListaGenratoriFormeInvarianti}
        E_{\delta},E_{-\delta},\omega^2_+,\omega^2_-,\omega^2,\omega^4_+,\omega^4_-,\omega^4\, 
    \end{equation}
and the dimensions of $\Lambda^i(\gm)^{\gh'}$ are 2, 4, 6, 9 and 12 for $i=1,2,3,4,5$, respectively.
\end{corollary}
\begin{proof}
    The dimensions are easily calculated by counting the occurrences of one--dimensional $\gh'$--modules and the computations can be conveniently carried out by means of the LiE program~\cite{LieProgram}.\par 
    By taking the wedge products of the basis elements obtained in   Proposition~\ref{propOmega1} and Proposition~\ref{propOmega2} we obtain all six  (linearly independent) generators of $\Lambda^3(\gm)^{\gh'}$, that are
    \begin{equation}
        E_{\delta} \wedge\omega^2_{\pm}\, ,\quad  E_{-\delta} \wedge\omega^2_{\pm}\, ,\quad E_{\pm \delta} \wedge\omega^2\, ,
    \end{equation}
    whereas in  $\Lambda^4(\gm)^{\gh'}$ we obtain six (linearly independent)  generators which, once the forms $\omega^4_\pm$ and $\omega^4$ obtained in Proposition~\ref{propFOURFORMS} are added, give all nine basis elements:
        \begin{eqnarray*}
        E_{\delta}\wedge E_{-\delta}\wedge\omega^2_{\pm}&=&E_{\delta}\wedge E_{-\delta}\wedge(E_{\pm\gamma_1}\wedge E_{\pm\gamma_2}-3E_{\pm\gamma_0}\wedge E_{\pm\gamma_3})\\
         E_{\delta}\wedge E_{-\delta}\wedge\omega^2&=&  E_{\delta}\wedge E_{-\delta}\wedge (3 E_{\gamma_0}\wedge E_{-\gamma_0}+E_{\gamma_1}\wedge E_{-\gamma_1}+E_{\gamma_2}\wedge E_{-\gamma_2}+3E_{\gamma_3}\wedge E_{-\gamma_3})\\ \omega^2_{+}\wedge\omega^2_{-}&=&
         E_{\gamma_1}\wedge E_{\gamma_2}\wedge E_{-\gamma_1}\wedge E_{-\gamma_2}-3(E_{\gamma_0}\wedge E_{\gamma_3}\wedge E_{-\gamma_1}\wedge    E_{-\gamma_2}+ E_{\gamma_1}\wedge E_{\gamma_2}\wedge E_{-\gamma_0}\wedge E_{-\gamma_3})\\
         &&+9 E_{\gamma_0}\wedge E_{\gamma_3}\wedge  E_{-\gamma_0}\wedge E_{-\gamma_3}
          \\
         \omega^2_{+}\wedge\omega^2&=& 3( E_{\gamma_1}\wedge E_{\gamma_2}\wedge  E_{\gamma_0}\wedge E_{-\gamma_0}+
         E_{\gamma_1}\wedge E_{\gamma_2}\wedge  E_{\gamma_3}\wedge E_{-\gamma_3}
        - E_{\gamma_0}\wedge E_{\gamma_3}\wedge  E_{\gamma_1}\wedge E_{-\gamma_1}\\
        &&-E_{\gamma_0}\wedge E_{\gamma_3}\wedge  E_{\gamma_2}\wedge E_{-\gamma_2}
 )  \\ \omega^2_{-}\wedge\omega^2&=&3( E_{-\gamma_1}\wedge E_{-\gamma_2}\wedge  E_{\gamma_0}\wedge E_{-\gamma_0}+
         E_{-\gamma_1}\wedge E_{-\gamma_2}\wedge  E_{\gamma_3}\wedge E_{-\gamma_3}\\
         &&- E_{-\gamma_0}\wedge E_{-\gamma_3}\wedge  E_{\gamma_1}\wedge E_{-\gamma_1}-E_{-\gamma_0}\wedge E_{-\gamma_3}\wedge  E_{\gamma_2}\wedge E_{-\gamma_2}
 ) \, ,\\
 \o^4_{\pm} &=& E_{\pm\gamma_0} \wedge E_{\pm\gamma_1}\wedge E_{\pm\gamma_2}\wedge E_{\pm \gamma_3 } \, ,\\
     \o^4 &=& E_{\gamma_0} \wedge E_{\gamma_1}\wedge E_{-\gamma_0}\wedge E_{-\gamma_1}+E_{\gamma_0} \wedge E_{\gamma_2}\wedge E_{-\gamma_0}\wedge E_{-\gamma_2}+\label{eq_MISSING_4form}\\
     &&+E_{\gamma_0} \wedge E_{\gamma_3}\wedge E_{-\gamma_1}\wedge E_{-\gamma_2}  +E_{\gamma_1} \wedge E_{\gamma_2}\wedge E_{-\gamma_0}\wedge E_{-\gamma_3} +\nonumber\\
     &&+E_{\gamma_1} \wedge E_{\gamma_3}\wedge E_{-\gamma_1}\wedge E_{-\gamma_3} +E_{\gamma_2} \wedge E_{\gamma_3}\wedge E_{-\gamma_2}\wedge E_{-\gamma_3}\, .
    \end{eqnarray*}
    The three four--forms $\omega^4_\pm$ and $\omega^4$ are exactly what is needed to fill up a basis. Finally, $\Lambda^5(\gm)^{\gh'}$ is generated by the  twelve basis elements
        \begin{equation}\label{eqLISTA_5FORME_INVARIANTI}
         \omega^2_{+}\wedge\omega^2_{-}\wedge E_{\pm \delta}\, ,\quad  \omega^2_{+}\wedge\omega^2\wedge E_{\pm \delta}\, ,\quad  \omega^2_{-}\wedge\omega^2\wedge E_{\pm \delta}
         \, ,\quad  \omega^4_{+}\wedge E_{\pm \delta}\, ,\quad  \omega^4_{-}\wedge E_{\pm \delta}\, ,\quad \omega^4\wedge E_{\pm \delta}\, .
    \end{equation}
\end{proof}


\bt  Forms~\eqref{eqListaGenratoriFormeInvarianti} 
define conformally invariant  forms  on  the bi-Lagrangian  symmetric manifold  $M^{10}_3 =  \mathsf{G}_2/\GL_2(\bR)^{\a_1}$.
The  pull--back  of these  forms   to the contact manifold   $ \Mcont_3 =\mathsf{G}_2/ \SL_2(\bR)^{\a_1}$
    generate   the  algebra  of    $\mathsf{G}_2$--invariant forms.
\et
\begin{proof}
    It suffices to observe that from 
\begin{eqnarray*}   \ad_{H_{\delta}}E_{\a_2 + k \a_1} &=&E_{\a_2 + k \a_1}\, ,\\  \ad_{H_{\delta}} E_{\delta}&=&2E_{\delta} \, ,
\end{eqnarray*}
it follows that
\begin{eqnarray*}   \ad_{H_{\delta}}E^*_{\pm\delta}&=& \mp 2 E_{\delta}^*\, ,\\ \ad_{H_{\delta}}\o_{\pm}^4  &=&  \mp 4 \o^4_{\pm}\, .
\end{eqnarray*}
\end{proof}

\begin{corollary}
 Any   $\mathsf{G}_2$--invariant  5--form  on the   contact manifold $\Mcont_3$ is a linear combination of the  forms~\eqref{eqLISTA_5FORME_INVARIANTI}.
\end{corollary} 

\subsection{Coordinate expressions of the $\mathsf{G}_2$--invariant MAEs}\label{secCoordTF}
Now we apply Theorem~\ref{thDarboux}, as well as the results recalled in Section~\ref{secMAEinCoordinates}, to find a coordinate expression of  the MAEs obtained above.  \par 
We define coordinates $\{x^0,x^1,x^2,x^3,x^4\}$ on the 5--dimensional homogeneous manifold
\begin{equation}
    F=G_2/P
\end{equation}
by means of the exponential map:
\begin{eqnarray}
    T_oF\simeq\gn_-&\longrightarrow& N^-\, ,\nonumber  \\
    X^-=x^0E_{-\gamma_0}+\cdots+x^3E_{-\gamma_3}+x^4E_{-\delta}&\longmapsto&\exp(X^-)\, .
\end{eqnarray}
\begin{remark}
    By switching the roles of $N^-$ and $N^+$ we obtain another chart: these two charts are enough to cover the whole $M$.
\end{remark}
Since $E_{\alpha}$ reads off the coefficient of $E_{-\alpha}$ in the above element $X^-$, we can make the following identifications
\begin{eqnarray*}
     dx^i&\longleftrightarrow& E_{\gamma_i}\, ,\quad i=0,1,2,3\, ,\\
     dx^4&\longleftrightarrow& E_{\delta}\, , \\
     du_i&\longleftrightarrow& E_{-\gamma_i}\, ,\quad i=0,1,2,3\, ,\\
     du_4&\longleftrightarrow& E_{-\delta}\, ,
\end{eqnarray*}
of covectors on the contact plane $\CC_o\equiv \gn^++\gn^-$ at the origin. With the above dictionary,  the twelve invariant 5--forms~\eqref{eqLISTA_5FORME_INVARIANTI}, that are covariant tensors on the vector space $\CC_o$, can be written down as degree--five (skew--symmetric) polynomials of the basis elements of $\CC_o^*$, see also the remark at the end of Section~\ref{secFiberWise}. These are, in that order:
\begin{eqnarray*}
    \omega^2_{+}\wedge\omega^2_{-}\wedge E_\delta&=&
         (dx^1\wedge dx^2\wedge du_1\wedge du_2-3(dx^0\wedge dx^3\wedge du_1\wedge    du_2+ dx^1\wedge dx^2\wedge du_0\wedge du_3)+9 dx^0\wedge dx^3\wedge  du_0\wedge du_3)\wedge {dx^4}\\
         \omega^2_{+}\wedge\omega^2_{-}\wedge E_{-\delta}&=&
         (dx^1\wedge dx^2\wedge du_1\wedge du_2-3(dx^0\wedge dx^3\wedge du_1\wedge    du_2+ dx^1\wedge dx^2\wedge du_0\wedge du_3)+9 dx^0\wedge dx^3\wedge  du_0\wedge du_3)\wedge {du_4}\\
          \omega^2_{+}\wedge\omega^2\wedge E_\delta&=& 3( dx^1\wedge dx^2\wedge  dx^0\wedge du_0+
         dx^1\wedge dx^2\wedge  dx^3\wedge du_3
        - dx^0\wedge dx^3\wedge  dx^1\wedge du_1-dx^0\wedge dx^3\wedge  dx^2\wedge du_2
 ) \wedge {dx^4}
 \\ 
          \omega^2_{+}\wedge\omega^2\wedge E_{-\delta}&=& 3( dx^1\wedge dx^2\wedge  dx^0\wedge du_0+
         dx^1\wedge dx^2\wedge  dx^3\wedge du_3
        - dx^0\wedge dx^3\wedge  dx^1\wedge du_1-dx^0\wedge dx^3\wedge  dx^2\wedge du_2
 ) \wedge {du_4}\\
 \omega^2_{-}\wedge\omega^2\wedge E_\delta&=&3( du_1\wedge du_2\wedge  dx^0\wedge du_0+
         du_1\wedge du_2\wedge  dx^3\wedge du_3
        - du_0\wedge du_3\wedge  dx^1\wedge du_1-du_0\wedge du_3\wedge  dx^2\wedge du_2
 )\wedge {dx^4} \, ,
 \\ \omega^2_{-}\wedge\omega^2\wedge E_{-\delta}&=&3( du_1\wedge du_2\wedge  dx^0\wedge du_0+
         du_1\wedge du_2\wedge  dx^3\wedge du_3
        - du_0\wedge du_3\wedge  dx^1\wedge du_1-du_0\wedge du_3\wedge  dx^2\wedge du_2
 )\wedge {du_4} \, ,\\
 \o^4_{+}\wedge E_\delta &=& dx^0 \wedge dx^1\wedge dx^2\wedge dx^3\wedge {dx^4} \, ,\\
  \o^4_{+}\wedge E_{-\delta} &=& dx^0 \wedge dx^1\wedge dx^2\wedge dx^3\wedge {du_4} \, ,\\
 \o^4_{-}\wedge E_\delta &=& du_0 \wedge du_1\wedge du_2\wedge du_3\wedge {dx^4} \, ,\\
 \o^4_{-}\wedge E_{-\delta} &=& du_0 \wedge du_1\wedge du_2\wedge du_3\wedge {du_4} \, ,\\
     \o^4 \wedge E_\delta&=&( dx^0 \wedge dx^1\wedge du_0\wedge du_1+dx^0 \wedge dx^2\wedge du_0\wedge du_2+\label{eq_MISSING_4form}\\
     &&+dx^0 \wedge dx^3\wedge du_1\wedge du_2  +dx^1 \wedge dx^2\wedge du_0\wedge du_3 +\nonumber\\
     &&+dx^1 \wedge dx^3\wedge du_1\wedge du_3 +dx^2 \wedge dx^3\wedge du_2\wedge du_3)\wedge {dx^4}\, ,\\
     \o^4\wedge E_{-\delta} &=& (dx^0 \wedge dx^1\wedge du_0\wedge du_1+dx^0 \wedge dx^2\wedge du_0\wedge du_2+\label{eq_MISSING_4form}\\
     &&+dx^0 \wedge dx^3\wedge du_1\wedge du_2  +dx^1 \wedge dx^2\wedge du_0\wedge du_3 +\nonumber\\
     &&+dx^1 \wedge dx^3\wedge du_1\wedge du_3 +dx^2 \wedge dx^3\wedge du_2\wedge du_3)\wedge {du_4}\, .
\end{eqnarray*}
Let us denote by  $\rho$  any of the twelve $(0,5)$--tensors on $\CC_o$ above: then there is a unique  $\mathsf{G}_2$--invariant five--form $\Omega^\rho$ on $M$, such that $\Omega^\rho_o=\rho$.  Then, formula~\eqref{eqDefMAE_general} allows to define twelve MAEs $\E_{\Omega^\rho}$, which can be now locally expressed in Darboux coordinates by means of the function $F_{\Omega^\rho}(u_{ij})$, cf.~\eqref{eqFunctionMAE}.\par
The global expression of $F_{\Omega^\rho}(u_{ij})$ can be obtained  from its restriction
\begin{equation}
    F_{\rho}:=F_{\Omega^\rho}|_{\LL(\CC_o,5)}
\end{equation}
to the fiber $\LL(\CC_o,5)$ of $N^{1}$ at $o\in N$: this  corresponds to the $\mathsf{G}_2$--equivariant extension
\begin{equation}
   \E_{\Omega^\rho}=\mathsf{G}_2\cdot \E_\rho:=\bigcup_{g\in\mathsf{G}_2}g(\E_\rho)\, ,\quad \E_\rho:=\{ F_{\rho}=0\}\, ,
\end{equation}
of the hypersurface $\E_\rho$ of $\LL(\CC_o,5)$.\par 
Below we compute the functions $F_\rho$ for all the twelve $(0,5)$--tensors above; to this end, 
we shall need the following symbols:
\begin{eqnarray*}
    M_{a}^l&\df&\textrm{rank--$4$ minor of $u_{ij}$ obtained by removing $a\Th$ row and $l\Th$ column,}\\
    M_{ab}^{lm}&\df&\textrm{rank--$3$ minor of $u_{ij}$ obtained by removing rows $a$, $b$,  and columns $l$, $m$,}\\
    M_{abc}^{lmn}&\df&\textrm{rank--$2$ minor of $u_{ij}$ obtained by removing rows $a$, $b$, $c$,  and columns $l$, $m$, $n$.}\\
\end{eqnarray*}
Then, the functions $F_\rho$ read  
\begin{eqnarray*}
\E_{\omega^2_{+}\wedge\omega^2_{-}\wedge E_\delta}&\longleftrightarrow& M_{124}^{034}-3\left(M_{034}^{034}+M_{124}^{124}\right)+9M^{124}_{034}=10M_{124}^{034}-3\left(M_{034}^{034}+M_{124}^{124}\right) \, ,\\
\E_{\omega^2_{+}\wedge\omega^2_{-}\wedge E_{-\delta}}&\longleftrightarrow& M_{12}^{03}-3\left(M_{03}^{03}+M_{12}^{12}\right)+9M^{12}_{03} = 10M_{12}^{03}-3\left(M_{03}^{03}+M_{12}^{12}\right)\, ,\\
\E_{\omega^2_{+}\wedge\omega^2 \wedge E_{\delta}}&\longleftrightarrow& 6(u_{03}+u_{12}) \, ,\\
\E_{\omega^2_{+}\wedge\omega^2 \wedge E_{-\delta}}&\longleftrightarrow& 6\left(M_{012}^{123}+M_{013}^{023}\right) \, ,\\
\E_{\omega^2_{-}\wedge\omega^2 \wedge E_{\delta}}&\longleftrightarrow& -3\left(M_{04}^{34}+M_{34}^{04}+M_{14}^{24}+M_{24}^{14}\right)=-6\left(M_{04}^{34} +M_{14}^{24} \right)  \, ,\\
\E_{\omega^2_{-}\wedge\omega^2 \wedge E_{-\delta}}&\longleftrightarrow& 3\left(M_0^3+M_3^0+M_1^2+M_2^1\right)=3\left(M_0^3 +M_1^2 \right)  \, ,\\
    \E_{\omega^4_+\wedge E_\delta}&\longleftrightarrow&\det(u_{ij})\, ,\\
    \E_{\omega^4_+\wedge E_{-\delta}}&\longleftrightarrow& u_{44}\, ,\\
    \E_{\omega^4_-\wedge E_\delta}&\longleftrightarrow& M_{4}^4 \, ,\\
    \E_{\omega^4_-\wedge E_{-\delta}}&\longleftrightarrow& \emptyset\, ,\\
    \E_{\omega^4 \wedge E_\delta}&\longleftrightarrow& M_{014}^{234}+M_{024}^{134}+M_{034}^{034}+M_{124}^{124}+M_{134}^{024}+M_{234}^{014} =2M_{014}^{234}+2M_{024}^{134}+M_{034}^{034}+M_{124}^{124}  \, ,\\
    \E_{\omega^4 \wedge E_{-\delta}}&\longleftrightarrow& M_{01}^{23}+M_{02}^{13}+M_{03}^{03}+M_{12}^{12}+M_{13}^{02}+M_{23}^{01}=2M_{01}^{23}+2M_{02}^{13}+M_{03}^{03}+M_{12}^{12}  \, .
\end{eqnarray*}
\subsection{Classification up to contactomorphisms}
We will need a particular element $\tau$ of the linear symplectic group $\Sp{\CC_o}=\Sp{\gm\oplus\gm^*}$: it is given by the  $10\times 10$ matrix
\begin{equation}
    \tau:=\left(\begin{matrix}
        0 & \id\\
        -\id & 0
    \end{matrix}\right)\, .
\end{equation}
Let us now observe that any element $g\in\mathsf{G}_2$ is, by construction, a contactomorphism of the manifold $N$, which preserves the integrable Lagrangian distributions $L^+$ and $L^-$ corresponding to the subalgebras $\gm$ and $\gm^*$. This means that the differential
\begin{equation}
d_og\in\Hom(\CC_o,\CC_{go})\simeq\Hom(\gm\oplus\gm^*,\gm\oplus\gm^*)
\end{equation}
preserves the symplectic structure of $\gm\oplus\gm^*$, preserves $\gm $ and $\gm^*$ separately and, moreover, the action on $\gm^*$ is dual to the action on $\gm$: this means that $d_og$ can be regarded as an element of $\Sp{\gm\oplus\gm^*}$, that is
\begin{equation}\label{eqDiGiO}
    d_og=\left(\begin{matrix}
        A & 0\\
        0 & (A^{-1})^t
    \end{matrix}\right)
\end{equation}
for some $5\times 5$ matrix $A$.\par
By employing the Darboux coordinates $\{x^0,x^1,x^2,x^3,x^4,u,u_0,u_1,u_2,u_3,u_4\}$ induced by the coordinates $\{x^0,x^1,x^2,x^3,x^4\}$ of $F$, see Section~\ref{secCoordTF}, we define the (total) Legendre transform
\begin{eqnarray*}
    \Phi:N &\longrightarrow & N\, ,\\
    (x^0,x^1,x^2,x^3,x^4,u,u_0,u_1,u_2,u_3,u_4)&\longmapsto &(u_0,u_1,u_2,u_3,u_4,u-x^iu_i,-x^0,-x^1,-x^2,-x^3,-x^4)\, ,
\end{eqnarray*}
which is a (local) contactomorphism of $N$, such that:
\begin{itemize}
    \item $\Phi$ sends the leaf $\mathcal{L}^\pm_p$ at $p\in N$ of the Lagrangian distribution $L^\pm$ to the leaf $\mathcal{L}^\mp_{\Phi(p)}$ at $\Phi(p)\in N$;
    \item the differential $d_p\Phi$, regarded as a linear  symplectomorphism of $\gm\oplus\gm^*$, coincides with $\tau$. 
\end{itemize}

\begin{theorem}\label{thTeoremaCONTROVERSO}
The MAE $\E_{\Omega^\rho}$ is contact--equivalent to the MAE $\E_{\Omega^{\tau(\rho)}}$.
\end{theorem}
\begin{proof}
It is enough to prove the identity
\begin{equation}\label{eqLegendrePutnuale}
    \mathsf{G}_2\cdot (\tau(\E_\rho))=\Phi(\mathsf{G}_2\cdot \E_\rho)\, .
\end{equation}
Since any point $p^1\in\E_\rho$ is interpreted as a Lagrangian plane $L_{p^1}$, the identity~\eqref{eqLegendrePutnuale} reads 
\begin{equation}\label{eqLegendrePutnualeBIS}
    \{ (d_og\circ\tau)(L_{p^1}) \mid  p^1\in\E_\rho\, , g\in \mathsf{G}_2  \}=\{(d_{\overline{g}o}\Phi\circ d_o\overline{g})(L_{p^1}) \mid  p^1\in\E_\rho\, , \overline{g}\in \mathsf{G}_2 \}\, .
\end{equation}
Let us observe that, for any $g\in G_2$, there exists another   $\overline{g}\in G_2$, such that 
\begin{equation}\label{eqDiGiBARO}
    d_o\overline{g}=\left(\begin{matrix}
        (A^{-1})^t & 0\\
        0 & A
    \end{matrix}\right) \, ,
    \end{equation}
Then, by using~\eqref{eqDiGiO} and~\eqref{eqDiGiBARO}, it is easy to show that
\begin{equation}\label{eqCommutativitaPezzotta}
    d_og\cdot\tau=\tau\cdot d_o\overline{g}\, .
\end{equation}
In light of the properties of $\Phi$, formula~\eqref{eqCommutativitaPezzotta} shows that~\eqref{eqLegendrePutnualeBIS} holds.\par
It only should be stressed that all the identifications we made in this proof are well defined up to an element of the stabilizer $\mathsf{SL}_2$ of $o$: this does not affect the final result, because the hypersurface $\E_\rho$ is $\mathsf{SL}_2$--invariant. 
\end{proof}
\begin{corollary}
A  $G_2$--invariant MAE that is obtained as the $G_2$--equivariant extension of  any of the   twelve hypersurfaces $\{F_\rho=0\}$ above is contactomorphic to the $G_2$--equivariant extension one of the following six $G_2$--invariant hypersurfaces:
\begin{eqnarray*}
 \textrm{quadratic \normalfont{(Q1)}:}&& 10M_{124}^{034}-3\left(M_{034}^{034}+M_{124}^{124}\right)=0\, ,\\
 \textrm{linear \normalfont{(L1)}:}&& u_{03}+u_{12}=0\, ,\\
 \textrm{quadratic \normalfont{(Q2)}:}&&M_{012}^{123}+M_{013}^{023}=0   \, ,\\
 \textrm{determinant \normalfont{(D)}:}&&\det(u_{ij})=0\, ,\\
 \textrm{linear \normalfont{(L2)}:}&& u_{44}=0\, ,\\
 \textrm{quadratic \normalfont{(Q3)}:}&&2M_{014}^{234}+2M_{024}^{134}+M_{034}^{034}+M_{124}^{124} =0  \, .
\end{eqnarray*}\end{corollary}
\begin{proof}
    Follows from the following identities, easily checked by direct computations:
\begin{eqnarray*}
    \tau(\E_{\omega^2_{+}\wedge\omega^2_{-}\wedge E_\delta})&=&\E_{\omega^2_{+}\wedge\omega^2_{-}\wedge E_{-\delta}}\, ,\\
    \tau( \E_{\omega^2_{+}\wedge\omega^2 \wedge E_{\delta}} ) &=& -\E_{\omega^2_{-}\wedge\omega^2 \wedge E_{-\delta}}  \, ,\\
     \tau(  \E_{\omega^2_{+}\wedge\omega^2 \wedge E_{-\delta}} )  &=& - \E_{\omega^2_{-}\wedge\omega^2 \wedge E_{\delta}}  \, ,\\
      \tau( \E_{\omega^4_+\wedge E_{\delta}} )   &=&\E_{\omega^4_-\wedge E_{-\delta}}   \, ,\\
       \tau( \E_{\omega^4_+\wedge E_{-\delta}}   )     &=& \E_{\omega^4_-\wedge E_{\delta}}  \, ,\\
         \tau(   \E_{\omega^4 \wedge E_\delta}  )     &=& \E_{\omega^4 \wedge E_{-\delta}}  \, .
\end{eqnarray*}
\end{proof}
\begin{theorem}
    The following results hold:
    \begin{enumerate}
        \item MAEs that are the $G_2$--equivariant extension of the hypersurfaces  labeled above (L1) and (Q2) are contact--equivalent;
        \item MAEs that are the $G_2$--equivariant extension of the hypersurfaces labeled above (D) and (L2) are contact--equivalent;
        \item MAEs that are the $G_2$--equivariant extension of the hypersurfaces labeled above (Q1) and (L1) are not contact--equivalent.
    \end{enumerate}
\end{theorem}
\begin{proof}
    Claims (1) and (2) can be proved analogously to the proof of Theorem~\ref{thTeoremaCONTROVERSO}, if $\tau$ is replaced by
    \begin{equation}
        \xi:=\left( \begin{matrix}
            \id_4 & 0 & 0& 0 \\
             0 & 0 & 0_4& 1 \\
              0_4 & 0 & \id_4& 0 \\
              0 & -1 & 0& 0 
        \end{matrix}   \right)\, ,
    \end{equation}
    and $\Phi$ is replaced by
    \begin{eqnarray*}
    \Xi:N &\longrightarrow & N\, ,\\
    (x^0,x^1,x^2,x^3,x^4,u,u_0,u_1,u_2,u_3,u_4)&\longmapsto &(x^0,x^1,x^2,x^3,u_4,u-x^4u_4,u_0,u_1,u_2,u_3,-x^4)\, .
\end{eqnarray*}
    Claim (3) can be proved by an analysis of the symbol. Since (L1) is linear, its symbol is constant and, as such, its rank never drops: it is constant to 4. On the other hand, the symbol of (Q1) is given, up to proportionality,  by the $5\times 5$ matrix
\begin{equation}
   \Smbl(\mathrm{Q1}) = \left(
\begin{array}{ccccc}
 -3 a_{12}^2 & 5 a_{10} a_{12} & -5 a_7 a_{12} & 3 a_3 a_{12} & 0 \\
 5 a_{10} a_{12} & -3 a_9 a_{12} & 3 a_6 a_{12} & -5 a_2 a_{12} & 0 \\
 -5 a_7 a_{12} & 3 a_6 a_{12} & -3 a_5 a_{12} & 5 a_1 a_{12} & 0 \\
 3 a_3 a_{12} & -5 a_2 a_{12} & 5 a_1 a_{12} & -3 a_3^2-3 a_6^2+10 a_2 a_7+3 a_5
   a_9-10 a_1 a_{10} & 0 \\
 0 & 0 & 0 & 0 & 0 \\
\end{array}
\right)\, ,
\end{equation}
where we have solved
\begin{equation}
    u_{00}=\frac{3 u_{{03}}^2+3 u_{{12}}^2-10 u_{{02}} u_{{13}}-3 u_{{11}}
   u_{{22}}+10 u_{{01}} u_{{23}}}{3 u_{{33}}}
\end{equation}
and we have set
\begin{equation}\label{eqCoordinateInterneQ1}
    a_1:=u_{01}\, ,\quad a_2:=u_{02}\, ,\ldots\, ,\quad  a_{14}:=u_{44}\, . 
\end{equation}
The rank of the above matrix is generically 4, but there are nonempty sets of points, where it drops.
\end{proof}


In terms of $u_{ij}$ variables, the list of $\mathsf{G_2}$--invariant MAEs, up to contact equivalence, reduces to   the $G_2$--equivariant extensions of the following hypersurfaces :
\begin{eqnarray*}
 \textrm{quadratic Q1:}&& 3 u_{03}^2+3 u_{12}^2-10 u_{02} u_{13}-3 u_{11}
   u_{22}+10 u_{01} u_{23}-3 u_{00} u_{33}=0\, ,\\
 \textrm{linear L1:}&& u_{03}+u_{12}=0\, ,\\
 \textrm{linear L2:}&& u_{44}=0\, ,\\
 \textrm{quadratic Q3:}&&2 u_{02}
   u_{13}+u_{11} u_{22}+2 u_{01} u_{23}+u_{00}
   u_{33}-u_{03}^2-4 u_{12} u_{03}-u_{12}^2=0  \, .
\end{eqnarray*}

\bibliographystyle{abbrvnat}
\bibliography{BibUniver}
\end{document}